\documentclass[10pt]{amsart}
\usepackage{amssymb}
\usepackage{amscd}

\theoremstyle{definition}
\newtheorem{thm}{Theorem}[section]
\newtheorem{lem}[thm]{Lemma}
\newtheorem{prp}[thm]{Proposition}
\newtheorem{dfn}[thm]{Definition}
\newtheorem{cor}[thm]{Corollary}

\newtheorem{ntn}[thm]{Notation}
\newtheorem{exa}[thm]{Example}

\newcommand{\beq}{\begin{equation}}
\newcommand{\eeq}{\end{equation}}
\newcommand{\beqr}{\begin{eqnarray*}}
\newcommand{\eeqr}{\end{eqnarray*}}
\newcommand{\bal}{\begin{align*}}
\newcommand{\eal}{\end{align*}}
\newcommand{\bei}{\begin{itemize}}
\newcommand{\eei}{\end{itemize}}
\newcommand{\limi}[1]{\lim_{{#1} \to \infty}}

\newcommand{\af}{\alpha}
\newcommand{\bt}{\beta}
\newcommand{\gm}{\gamma}
\newcommand{\dt}{\delta}
\newcommand{\ep}{\varepsilon}

\newcommand{\et}{\eta}

\newcommand{\io}{\iota}
\newcommand{\te}{\theta}
\newcommand{\ld}{\lambda}
\newcommand{\sm}{\sigma}

\newcommand{\ph}{\varphi}
\newcommand{\ps}{\psi}
\newcommand{\rh}{\rho}

\newcommand{\ta}{\tau}

\newcommand{\Gm}{\Gamma}
\newcommand{\Dt}{\Delta}

\newcommand{\Ld}{\Lambda}

\newcommand{\Z}{{\mathbf{Z}}}
\newcommand{\R}{{\mathbf{R}}}
\newcommand{\C}{{\mathbf{C}}}
\newcommand{\N}{{\mathbf{N}}}

\newcommand{\id}{{\mathrm{id}}}

\newcommand{\dist}{{\mathrm{dist}}}
\newcommand{\sa}{{\mathrm{sa}}}
\newcommand{\spec}{{\mathrm{sp}}}
\newcommand{\Ad}{{\mathrm{Ad}}}
\newcommand{\Aff}{{\mathrm{Aff}}}
\newcommand{\Aut}{{\mathrm{Aut}}}
\newcommand{\diag}{{\mathrm{diag}}}
\newcommand{\supp}{{\mathrm{supp}}}

\newcommand{\card}{{\mathrm{card}}}

\newcommand{\Mi}{M_{\infty}}

\newcommand{\andeqn}{\,\,\,\,\,\, {\mbox{and}} \,\,\,\,\,\,}

\newcommand{\ts}[1]{{\textstyle{#1}}}
\newcommand{\ds}[1]{{\displaystyle{#1}}}
\newcommand{\ssum}[2]{{\ts{ {\ds{\sum}}_{#1}^{#2} }}}
\newcommand{\susum}[1]{{\ts{ {\ds{\sum}}_{#1} }}}
\newcommand{\scup}[2]{{\ts{ {\ds{\bigcup}}_{#1}^{#2} }}}
\newcommand{\sucup}[1]{{\ts{ {\ds{\bigcup}}_{#1} }}}

\newcommand{\ca}{C*-algebra}
\newcommand{\ct}{continuous}
\newcommand{\pj}{projection}

\newcommand{\hm}{homomorphism}
\newcommand{\fd}{finite dimensional}
\newcommand{\wolog}{without loss of generality}
\newcommand{\Wolog}{Without loss of generality}

\newcommand{\ifo}{if and only if}

\newcommand{\mops}{mutually orthogonal \pj s}

\newcommand{\hme}{homeomorphism}

\newcommand{\cfn}{continuous function}
\newcommand{\hsa}{hereditary subalgebra}

\newcommand{\mvnt}{Murray-von Neumann equivalent}

\newcommand{\tRp}{tracial Rokhlin property}

\newcommand{\sfsuca}{stably finite simple unital \ca}

\newcommand{\rsz}[1]{\raisebox{0ex}[0.8ex][0.8ex]{$#1$}}

\title[Rank of crossed products]{Stable and real rank for
    crossed products by automorphisms with the tracial Rokhlin
    property}

\author{Hiroyuki Osaka}
\author{N.\  Christopher Phillips}

\date{30 August 2004}

\address{Department of Mathematics, Ritsumeikan University,
       Kusatsu, Shiga, 525-8577, Japan.}
\email[]{osaka@se.ritsumei.ac.jp}

\address{Department of Mathematics, University  of Oregon,
       Eugene OR 97403-1222, USA.}
\email[]{ncp@darkwing.uoregon.edu}

\subjclass[2000]{Primary 46L55;
   Secondary 16S35, 19A13, 46L35, 46L40.}
\thanks{Research partially supported by
JSPS Grant for Scientific Research No.\  14540217(c)(1),
and by NSF grants DMS-0070776 and DMS-0302401.}

\begin{document}

\begin{abstract}
We introduce the \tRp\  for automorphisms of
stably finite simple unital \ca s containing enough \pj s.
This property is formally weaker than the various Rokhlin properties
considered by Herman and Ocneanu, Kishimoto, and Izumi.
Our main results are as follows.
Let $A$ be a stably finite simple unital \ca,
and let $\af$ be an automorphism of $A$ which has the \tRp.
Suppose $A$ has real rank zero and stable rank one,
and suppose that the order on \pj s over $A$ is determined
by traces.
Then the crossed product algebra
$C^* (\Z, A, \af)$ also has these three properties.

We also present examples of \ca s $A$ with automorphisms $\af$
which satisfy the above assumptions, but such that
$C^* (\Z, A, \af)$ does not have tracial rank zero.
\end{abstract}

\maketitle

\setcounter{section}{-1}

\section{Introduction}

\indent
We introduce the \tRp\  for automorphisms of
stably finite simple unital \ca s containing enough \pj s.
This property is formally weaker than the various
Rokhlin properties which have appeared in the literature,
such as in~\cite{HO}, \cite{Ks3}, and~\cite{Iz},
at least for \ca s which are tracially~AF in the sense of~\cite{LnTAF},
in roughly the same way that being tracially~AF is weaker
than the local characterization of AF~algebras
(Theorem~2.2 of~\cite{Brt}).

Our main results are as follows.
Let $A$ be a stably finite simple unital \ca,
and let $\af$ be an automorphism of $A$ which has the \tRp.
Suppose $A$ has real rank zero and stable rank one,
and suppose that the order on \pj s over $A$ is determined
by traces (Blackadar's Second Fundamental
Comparability Question, 1.3.1 of~\cite{Bl3}, for $M_{\infty} (A)$).
Then $C^* (\Z, A, \af)$ also has these three properties.
In fact, we will see that not all the hypotheses on $A$ are
needed for all the conclusions.

The proofs are adapted from~\cite{Ph10}.
The arguments here are more difficult for several reasons.
First, in~\cite{Ph10} there is a single ``large'' AF~subalgebra,
and the properties of the reduced groupoid \ca\  are obtained
by comparison with this subalgebra.
In this paper, we are not able to choose nested approximating
subalgebras;
moreover, even if we were, the direct limit would not be~AF.
Second,
we assume that the order on \pj s over $A$
is determined by all traces,
but only the invariant traces extend to the crossed product.
Third, in~\cite{Ph10} we relied on previous work
to get from the Rokhlin property to the existence of appropriate
subalgebras, but in the present paper we must do the analogous
construction from scratch.

Kishimoto has proved (Theorem~6.4 of~\cite{Ks4})
that if $A$ is a simple unital AT~algebra with real rank zero
which has a unique tracial state,
and if $\af \in \Aut (A)$ has the Rokhlin property
and satisfies a kind of approximate innerness assumption,
then $C^* (\Z, A, \af)$ is again
a simple unital AT~algebra with real rank zero.
As this paper was in progress, H.\  Lin and the first
author have generalized this~\cite{LO}.
Let $A$ be a simple separable unital C*-algebra
which satisfies the Universal Coefficient Theorem,
which has tracial rank zero,
and which has a unique tracial state.
If $\af \in \Aut (A)$ has the Rokhlin property
and if $\af^n$ is an approximately inner for some $n > 0,$
then $C^* (\Z, A, \af)$ is a simple AH~algebras with no
dimensional growth and real rank zero.
It seems reasonable to hope that whenever $A$ is tracially~AF
and $\af$ has the \tRp,
then $C^* (\Z, A, \af)$ is again tracially~AF.
However, the theorems in this paper also apply to
automorphisms of \ca s which are not tracially~AF,
and for which the crossed products are also not tracially~AF.
We give some examples in Section~\ref{Sec:NonFurstEx}.

Our motivating examples are the noncommutative Furstenberg
transformations, which are automorphisms of the irrational rotation
algebras analogous to Furstenberg transformations on the torus,
and the automorphisms in the crossed product descriptions
of the simple quotients of the \ca s of discrete subgroups
of nilpotent Lie groups studied in~\cite{MW2} and~\cite{MW4}.
These automorphisms do not satisfy the hypotheses in~\cite{LO},
although in these cases we believe that
the crossed products are in fact tracially~AF.
We treat these examples in a separate paper~\cite{OP2}.
In that paper we also show that an automorphism of a
simple unital tracially~AF \ca\ $A$  with unique trace $\ta$
has the \tRp\  \ifo\  all nontrivial powers of
the corresponding automorphism of the factor $\pi_{\ta} (A)'',$
obtained from the
Gelfand-Naimark-Segal representation associated with $\ta,$
are outer.

This paper is organized as follows.
In Section~\ref{Sec:tRp} we introduce the \tRp,
and prove, under reasonable conditions,
that it is implied by various forms of the Rokhlin property
in the literature.
We also obtain several elementary consequences.
In Section~\ref{Sec:Tower},
we use the \tRp\  to construct ``large'' subalgebras of
$C^* (\Z, A, \af)$ which are stably isomorphic to $A.$
The next three sections treat, in order,
the order on \pj s, real rank zero, and stable rank one.
These are the analogs of Sections~3, 4, and~5 of~\cite{Ph10}.
It is in Section~\ref{Sec:RR} that the weaker conditions
satisfied by the subalgebras cause the greatest additional
difficulty.
We also obtain several other results:
the restriction map from tracial states on $C^* (\Z, A, \af)$
to invariant tracial states on $A$ is bijective,
and, under suitable hypotheses,
$C^* (\Z, A, \af)$ satisfies the local approximation property
of Popa~\cite{Pp}.
The last section gives some examples,
but the ones we are most interested in are in~\cite{OP2}.

The first author is grateful to Masaki Izumi
for valuable discussions,
and second author is grateful to Nate Brown, Masaki Izumi,
Cornel Pasnicu,
Christian Skau, and Takeshi Katsura for valuable discussions.

\section{The tracial Rokhlin property}\label{Sec:tRp}

\indent
We begin by defining the \tRp\  for single automorphisms
(actions of $\Z$).
It is closely related to, but slightly weaker than,
the approximate Rokhlin property of Definition~4.2 of~\cite{Ks1}.
To our knowledge, the idea was first introduced in~\cite{BEK}.
It is closely related to the \tRp\  for actions of finite
cyclic groups, as in~\cite{Ph11}.

\begin{dfn}\label{TRPDfn}
Let $A$ be a \sfsuca\   and let $\af \in \Aut (A).$
We say that $\af$ has the {\emph{tracial Rokhlin property}}
if for every finite set $F \subset A,$ every $\ep > 0,$
every $n \in \N,$
and every nonzero positive element $x \in A,$
there are \mops\  $e_0, e_1, \ldots, e_n \in A$ such that:
\begin{itemize}
\item[(1)]
$\| \af (e_j) - e_{j + 1} \| < \ep$ for $0 \leq j \leq n - 1.$
\item[(2)]
$\| e_j a - a e_j \| < \ep$ for $0 \leq j \leq n$ and all $a \in F.$
\item[(3)]
With $e = \sum_{j = 0}^{n} e_j,$ the \pj\  $1 - e$ is \mvnt\  to a
\pj\  in the \hsa\  of $A$ generated by $x.$
\end{itemize}
\end{dfn}

We do not say anything about $\af (e_n).$

It is not completely clear that Condition~(3) is the right
condition in the general case.
We return to this point,
and to the comparison of our definition with others,
after some preliminaries.

\begin{ntn}\label{TraceNtn}
Let $A$ be a unital \ca.
We denote by $T (A)$ the set of all tracial states on $A,$
equipped with the weak* topology.
For any element of $T (A),$
we use the same letter for its standard extension to $M_n (A)$
for arbitrary $n,$
and to $\Mi (A) = \bigcup_{n = 1}^{\infty} M_n (A)$ (no closure).
\end{ntn}

\begin{dfn}\label{OrdDetD}
Let $A$ be a unital \ca.
We say that the {\emph{order on \pj s over $A$ is determined by traces}}
if whenever $n \in \N$ and $p, q \in \Mi (A)$ are \pj s such that
$\ta (p) < \ta (q)$ for all $\ta \in T (A),$
then $p \precsim q.$
\end{dfn}

This is Blackadar's Second Fundamental Comparability Question
for $\Mi (A).$
See 1.3.1 in~\cite{Bl3}.

In all applications so far,
in addition to the conditions in Definition~\ref{TRPDfn},
the algebra $A$ has real rank zero (Section~1 of~\cite{BP}),
and the order on \pj s over $A$ is determined by traces.
In this case, we can replace the third condition by one involving
traces:

\begin{lem}\label{DfnUsingTrace}
Let $A$ be a \sfsuca\   such that ${\mathrm{RR}} (A) = 0$
and the order on \pj s over $A$ is determined by traces.
Let $\af \in \Aut (A).$
Then $\af$ has the \tRp\  %
\ifo\  for every finite set $F \subset A,$ every $\ep > 0,$
and every $n \in \N,$
there are \mops\  $e_0, e_1, \ldots, e_n \in A$ such that:
\begin{itemize}
\item[(1)]
$\| \af (e_j) - e_{j + 1} \| < \ep$ for $0 \leq j \leq n - 1.$
\item[(2)]
$\| e_j a - a e_j \| < \ep$ for $0 \leq j \leq n$ and all $a \in F.$
\item[(3)]
With $e = \sum_{j = 0}^{n} e_j,$ we have $\ta (1 - e) < \ep$
for all $\ta \in T (A).$
\end{itemize}
\end{lem}

\begin{proof}
If the condition of Definition~\ref{TRPDfn} holds
and $\ep,$ $n,$ and $F$ are given,
then we can use Theorem~1.1(a) of~\cite{Zh7} to find a \pj\  $x \in A$
such that $\ta (x) < \ep$ for all $\ta \in T (A).$
Then apply Definition~\ref{TRPDfn} with this $x$
and with $\ep,$ $n,$ and $F$ as given.
Conversely, assume the condition of the lemma,
and let $\ep,$ $n,$ $F,$ and $x$ be given.
Choose a nonzero \pj\  $q \in {\overline{x A x}},$
and apply the condition of the lemma with $\ep$ replaced by
$\min \left( \ep, \, \inf_{\ta \in T (A)} \ta (q) \right).$
The assumption that the order on \pj s over $A$ is determined by traces
implies that $1 - e \precsim q,$ giving~(3) of Definition~\ref{TRPDfn}.
\end{proof}

We now want to relate the \tRp\  to forms of the Rokhlin property
which have appeared in the literature.
The most important of these is as follows.
(See, for example, Definition~2.5 of~\cite{Iz},
and Condition~(3) in Proposition~1.1 of~\cite{Ks3}.)

\begin{dfn}\label{StrictRP}
Let $A$ be a simple unital \ca\   and let $\af \in \Aut (A).$
We say that $\af$ has the {\emph{Rokhlin property}}
if for every finite set $F \subset A,$ every $\ep > 0,$
every $n \in \N,$
there are \mops\  %
\[
e_0, \, e_1, \, \ldots, \, e_{n - 1},
       \, f_0, \, f_1, \, \ldots, \, f_n \in A
\]
such that:
\begin{itemize}
\item[(1)]
$\| \af (e_j) - e_{j + 1} \| < \ep$ for $0 \leq j \leq n - 2$
and $\| \af (f_j) - f_{j + 1} \| < \ep$ for $0 \leq j \leq n - 1.$
\item[(2)]
$\| e_j a - a e_j \| < \ep$ for $0 \leq j \leq n - 1$ and all $a \in F,$
and
$\| f_j a - a f_j \| < \ep$ for $0 \leq j \leq n$ and all $a \in F.$
\item[(3)]
$\sum_{j = 0}^{n - 1} e_j + \sum_{j = 0}^{n} f_j = 1.$
\end{itemize}
\end{dfn}

We will also consider analogs of the original version for \ca s,
in for example the first definition of~\cite{HO},
in which Condition~(3) of Definition~\ref{TRPDfn}
is replaced by $\sum_{j = 0}^{n} e_j = 1,$
but in which, as in Lemma~\ref{FOrdOnTr} below,
the towers are only required to exist for all $n$ in
an unbounded subset $S \subset \N$ which does not depend on $\ep$ and $F.$

We do not know whether these properties imply the \tRp\  in
the generality we have considered so far.
We prove that they do under the following sets of hypotheses,
in all of which we assume that $A$ is stably finite, simple, and unital:
\begin{itemize}
\item
${\mathrm{RR}} (A) = 0,$
the order on \pj s over $A$ is determined by traces,
and the \hme\  $\ta \mapsto \ta \circ \af$ of $T (A)$
has finite order.
\item
$A$ is approximately divisible in the sense of~\cite{BKR},
all quasitraces on $A$ are traces,
and \pj s in $A$ distinguish the traces on $A.$
\item
$A$ has tracial rank zero.
\end{itemize}
Together, these cover most of the interesting cases in which
$A$ has real rank zero.
Note that $\ta \mapsto \ta \circ \af$ has finite order
whenever all tracial states are $\af$-invariant
(in particular, whenever $\af$ is approximately inner or
${\mathrm{RR}} (A) = 0$ and $\af$ is trivial on K-theory),
and also whenever there are only finitely many extreme tracial states.

We can obtain a version of the \tRp\  which is implied by
the Rokhlin property in full generality
by allowing two Rokhlin towers in Definition~\ref{TRPDfn},
as is done in Definition~\ref{StrictRP},
but still allowing a remainder \pj.
The proofs of our main theorems
should all still work under this condition.
Another possibility, motivated by Proposition~\ref{OrdPjCrPrd},
is to merely require that there be $q \in {\overline{x A x}}$
such that $1 - e \sim q,$
with equivalence in $C^* (\Z, A, \af)$ rather than $A.$
We have not checked whether our proofs still work with this
assumption.
Our motivation for using the definition as stated
is Theorem~2.14 of~\cite{OP2},
which under certain conditions relates the \tRp\  to
a property having the form of the Rokhlin property for
automorphisms of factors of type~II$_1.$

There are immediate K-theoretic obstructions to any version of
the Rokhlin property involving only one tower and
requiring $\sum_{j = 0}^{n} e_j = 1,$ as in~\cite{HO}.
However, we know of no K-theoretic obstructions for an
automorphism of a simple unital \ca\  with real rank zero to
have the Rokhlin property as in Definition~\ref{StrictRP}.
It is in fact implicit in several proofs in the literature,
in particular in the proof of Theorem~4.1 of~\cite{Ks1},
that, under suitably restrictive hypotheses,
the tracial Rokhlin property implies the Rokhlin property.
The hypotheses include the assumption that
there are no infinitesimals in the $K_0$-group.
We know of no examples of automorphisms with the \tRp\  which do
not have the Rokhlin property of Definition~\ref{StrictRP}.

Kishimoto's definition of the approximate Rokhlin property,
Definition~4.2 of~\cite{Ks1}, specifies that $A$ is an AF algebra,
and, instead of having a finite set $F,$ it assumes a \fd\  subalgebra
$B$ is given; in place of approximate commutativity with a finite set,
it requires
that every $e_j$ commute exactly with every element of $B.$
More significantly, that definition also requires that
$\| \af (e_n) - e_0 \|$ be small,
while we make no assumption on $\af (e_n).$
It furthermore has a slightly different version of Condition~(3).
When the order on \pj s is determined by traces,
the analog of the approximate Rokhlin property in our case
is formally stronger than the \tRp\  as we defined it.
In Lemma~4.4 of~\cite{Ks1}, Kishimoto explicitly proves
that on a simple unital AF algebra whose $K_0$-group is finitely
generated and contains no infinitesimal elements,
the approximate Rokhlin property implies the Rokhlin property.

We now prove that the Rokhlin property implies the \tRp\  under
the first of the sets of hypotheses discussed above.
We need a lemma.

\begin{lem}\label{FOrdOnTr}
Let $A$ be a \sfsuca\   such that ${\mathrm{RR}} (A) = 0$
and the order on \pj s over $A$ is determined by traces.
Let $\af \in \Aut (A),$
and suppose that the \hme\  of $T (A)$
given by $\ta \mapsto \ta \circ \af$ has finite order.
Then $\af$ has the \tRp\  %
\ifo\  there is an unbounded set $S \subset \N$ and a constant $C > 0$
such that for every finite set $F \subset A,$ every $\ep > 0,$
and every $n \in S,$
there are \mops\  $e_0, e_1, \ldots, e_n \in A$ such that:
\begin{itemize}
\item[(1)]
$\| \af (e_j) - e_{j + 1} \| < \ep$ for $0 \leq j \leq n - 1.$
\item[(2)]
$\| e_j a - a e_j \| < \ep$ for $0 \leq j \leq n$ and all $a \in F.$
\item[(3)]
With $e = \sum_{j = 0}^{n} e_j,$
we have $\ta (1 - e) \leq C (n + 1)^{- 1}$ for all $\ta \in T (A).$
\end{itemize}
\end{lem}

\begin{proof}
That the \tRp\  implies the condition of the lemma is clear from
Lemma~\ref{DfnUsingTrace}.
Conversely, assume the conditions of the lemma.
We prove the condition in Lemma~\ref{DfnUsingTrace}
for $\ep,$ $n,$ and $F$ as there.
\Wolog\  $\ep < 1.$
Let $k$ be a positive integer such that
$\ta \circ \af^k = \ta$ for all $\ta \in T (A).$
Thus, if $p_0, p_1, \ldots, p_N$ are \pj s such that
$\| \af (p_j) - p_{j + 1} \| < 1$ for $l \leq j \leq l + k - 1,$
then $\ta (p_{l + k}) = \ta (p_l)$ for all $\ta \in T (A).$
Choose $N \in S$ with
\[
N \geq
  \max \left( \frac{2 k (n + 1)}{\ep}, \, \frac{2 C}{\ep} \right).
\]
Apply the condition of this lemma with $N$ in place of $n$ and with
$\ts{ \frac{1}{2}} \ep N^{- 1}$ in place of $\ep,$
to find \mops\  $p_0, p_1, \ldots, p_N \in A,$
and set $p = \sum_{m = 0}^{N} p_m.$
Thus $\ta (1 - p) \leq C (N + 1)^{- 1} < \ts{ \frac{1}{2}} \ep.$
Write $N + 1 = r k (n + 1) + s$
with $r \in \N$ and $0 \leq s < k (n + 1),$
so that $r > 2 \ep^{-1}.$
For $0 \leq j \leq n$ set
\[
e_j = p_j + p_{j + n + 1} + \cdots + p_{j + (r k - 1) (n + 1)},
\]
and set $e = \sum_{j = 0}^{n} e_j.$
We easily get
\[
\| \af (e_j) - e_{j + 1} \| < r k \cdot \ts{ \frac{1}{2}} \ep N^{-1}
  \leq \ep
\andeqn
\| e_j a - a e_j \| < r k \cdot \ts{ \frac{1}{2}} \ep N^{-1}
  \leq \ep
\]
for $0 \leq j \leq n$ and $a \in F.$
For the trace estimate, for $0 \leq l \leq r - 1$ set
\[
q_l = p_{k (n + 1) l} + p_{k (n + 1) l + 1}
              + \cdots + p_{k (n + 1) (l + 1) - 1}.
\]
Then
\[
e = q_0 + q_1 + \cdots + q_{r - 1}
\andeqn
p - e = p_{k (n + 1) r} + p_{k (n + 1) r + 1}
             + \cdots + p_N.
\]
Using periodicity of $\af$ on $T (A)$ of order $k (n + 1),$
and $N + 1 - k (n + 1) r < k (n + 1),$
we get for every $\ta \in T (A)$
\[
\ta (q_0) = \ta (q_1) = \cdots = \ta (q_{r - 1}) > \ta (p - e),
\]
so
\[
\ta (p - e) < \frac{1}{r + 1} < \frac{\ep}{2}.
\]
Thus $\ta (1 - e) < \ep.$
\end{proof}

\begin{prp}\label{SRPImpTRP}
Let $A$ be a \sfsuca\   such that ${\mathrm{RR}} (A) = 0$
and the order on \pj s over $A$ is determined by traces.
Let $\af \in \Aut (A),$
and suppose that the \hme\  of $T (A)$
given by $\ta \mapsto \ta \circ \af$ has finite order.
If $\af$ has the Rokhlin property of Definition~\ref{StrictRP},
then $\af$ has the \tRp.
\end{prp}

\begin{proof}
Let $k$ be a positive integer such that
$\ta \circ \af^k = \ta$ for all $\ta \in T (A).$
We verify the condition of Lemma~\ref{FOrdOnTr},
with
\[
S = \{ k - 1, \, 2 k - 1, \, \ldots \} \andeqn C = k.
\]
Let $r \in \N.$
Apply the Rokhlin property
with $\min \left( 1, \ts{\frac{1}{2}} \ep \right)$ in place of $\ep,$
with $F$ as given, and with $n = r k.$
Let
\[
p_0, \, p_1, \, \ldots, \, p_{r k - 1},
       \, q_0, \, q_1, \, \ldots, \, q_{r k} \in A
\]
be the resulting \pj s,
and take $e_j = p_j + q_j$ for $0 \leq j \leq r k - 1.$
Then $1 - \sum_{j = 0}^{r k - 1} e_j = q_{r k}.$
By periodicity of $\af$ on $T (A),$ we have, for all $\ta \in T (A),$
\[
\ta (q_{r k}) = \ta (q_0)
 \leq \ta (q_0 + q_1 + \cdots + q_{k - 1})
 < \frac{1}{r} = \frac{C}{r k}.
\]
This completes the proof.
\end{proof}

Now we prove that the Rokhlin property implies the \tRp\  under
the second and third of the sets of hypotheses discussed above.
In the proofs above, the ``leftover \pj'' in the \tRp\  was
the sum of the \pj s in a small part of the tower obtained
from the Rokhlin property.
Without something like $\ta \mapsto \ta \circ \af$ having finite
order, we don't see how to make such a proof work.
Instead, we must divide a tower in parallel towers of the
same height, and arrange to omit different \pj s in each,
so that altogether the ``leftover'' consists of a small part
of each of the \pj s in the original towers.
This is a bit messy to write down.

We need a preparatory lemma for each set of hypotheses.

\begin{lem}\label{ZhangTAF}
Let $A$ be a simple unital infinite dimensional \ca\  %
with tracial rank zero.
Let $p \in A$ be a nonzero \pj, let $F \subset p A p$ be a finite set,
let $m \in \N$ be a power of two, and let $\ep > 0.$
Then there exist \pj s $p_0, p_1, \ldots, p_m \in A$ such that
\[
\sum_{r = 0}^m p_r = p, \,\,\,\,\,\,
p_1 \sim p_2 \sim \cdots \sim p_m,
\andeqn
p_0 \precsim p_1,
\]
and such that $\| [p_r, a] \| < \ep$
for $0 \leq r \leq m$ and all $a \in F.$
\end{lem}

\begin{proof}
We have ${\mathrm{RR}} (A) = 0$ by Theorem~3.4 of~\cite{LnTAF}.
It now follows easily from Theorem~1.1(a) of~\cite{Zh7}
that there is a nonzero \pj\  $q \in p A p$
such that $p A p$ contains $2 m + 3$ \mops, each \mvnt\  to $q.$
It follows from Theorem~3.12 of~\cite{LnTAF} that $p A p$
also has tracial rank zero.
Therefore there is a finite dimensional subalgebra $E \subset p A p,$
with identity $e \leq p,$ such that
$\| [e, a] \| < \frac{1}{6} \ep$ for all $a \in F,$
such that for every $a \in F$ there is $b \in E$
with $\| b - e a e \| < \frac{1}{6} \ep,$
and such that $p - e \precsim q.$

Let
\[
B = E' \cap e A e
 = \{ x \in e A e \colon {\mbox{$x c = c x$ for every $c \in E$}} \}.
\]
Write $e = \sum_{k = 1}^n e_k$
as a sum of minimal central \pj s of $E.$
Let $f_k \in E$ be a minimal \pj\  with $f_k \leq e_k.$
Then $B = \bigoplus_{k = 1}^n e_k B e_k,$
and $e_k B e_k \cong f_k A f_k$ is simple and has real rank zero.
Since $2 m$ is also a power of two,
Theorem~1.1(a) of~\cite{Zh7} therefore provides
\pj s $q_{k, 0}, q_{k, 1}, \ldots, q_{k, 2 m} \in e_k B e_k$ such that
\[
\sum_{r = 0}^{2 m} q_{k, r} = e_k, \,\,\,\,\,\,
q_{k, 1} \sim q_{k, 2} \sim \cdots \sim q_{k, 2 m},
\andeqn
q_{k, 0} \precsim q_{k, 1}.
\]
Then define
\[
p_0 = p - e + \sum_{k = 1}^n q_{k, 0},
\]
and, for $0 \leq r \leq m,$
\[
p_r = \sum_{k = 1}^n (q_{k, \, 2 r - 1} + q_{k, \, 2 r}).
\]
We prove that these \pj s satisfy the conclusion of the lemma.

It is clear that
$\sum_{r = 0}^m p_r = p$ and $p_1 \sim p_2 \sim \cdots \sim p_m.$
To prove that $p_0 \precsim p_1,$
we use the fact that,
by Theorems~5.8 and~6.8 of~\cite{LnTTR},
the order on projections over $A$ is determined by traces.
We certainly have
$\sum_{k = 1}^n q_{k, 0} \precsim \sum_{k = 1}^n q_{k, 1}.$
Let $\ta \in T (A).$
Set $\bt = \sum_{k = 1}^n \ta (q_{k, 2}).$
Then $\sum_{k = 1}^n \ta (q_{k, 0}) \leq \bt$
and $\ta (p_r) = 2 \bt$ for $1 \leq r \leq m.$
So
\[
\ta (p) = \ta (p - e) + \sum_{k = 1}^n \ta (q_{k, 0}) + 2 m \bt
   \leq \ta (p - e) + (2 m + 1) \bt.
\]
Using the choice of $q,$ we get
\[
\ta (p - e) \leq \ta (q) < \frac{\ta (p)}{2 m + 2}.
\]
It follows that $\bt > \ta (p) / (2 m + 2).$
Therefore $\ta (p - e) < \bt.$
Since this is true for all $\ta \in T (A),$
we conclude that $p - e \precsim \sum_{k = 1}^n q_{k, 2}.$
Combining that with our first observation gives $p_0 \precsim p_1,$
as desired.

It remains to estimate $\| [p_r, a] \|$
for $0 \leq r \leq m$ and $a \in F.$
Choose $b \in E$ such that $\| b - e a e \| < \frac{1}{6} \ep.$
Then
\[
\| [ b + (p - e) a (p - e)] - a \|
  \leq \| e a e - b \| + \| (p - e) a e \| + \| e a (p - e) \|
  < 3 \left( \tfrac{1}{6} \ep \right)
  = \tfrac{1}{2} \ep,
\]
and $b + (p - e) a (p - e)$ commutes with $p_r,$
so $\| [p_r, a] \| < \ep.$
\end{proof}

\begin{lem}\label{ZhangAppDiv}
Let $A$ be a simple separable unital approximately divisible \ca.
Let $p \in A$ be a nonzero \pj, let $F \subset p A p$ be a finite set,
let $m \in \N,$ and let $\ep > 0.$
Then there exist \pj s $p_0, p_1, \ldots, p_m \in A$ such that
\[
\sum_{r = 0}^m p_r = p, \,\,\,\,\,\,
p_1 \sim p_2 \sim \cdots \sim p_m,
\andeqn
p_0 \precsim p_1,
\]
and such that $\| [p_r, a] \| < \ep$
for $0 \leq r \leq m$ and all $a \in F.$
\end{lem}

\begin{proof}
It follows from Corollary~2.10 of~\cite{BKR}
that there is a finite dimensional unital subalgebra
$E \subset A$ such that $\| [a, x] \| \leq \frac{1}{2} \ep \| x \|$
for all $x \in E$ and $a \in F,$
and such that $E \cong \bigoplus_{l = 1}^t M_{n (l)}$
with $n (l) \geq m^2$ for $1 \leq l \leq t.$
Let $( e_{j, k}^{(l)} )_{1 \leq j, k \leq n (l)}$
be a system of matrix units for the $l$-th summand of $E.$
Write $n (l) = d (l) m + z (l)$ with $0 \leq z (l) \leq m - 1.$
Note that $d (l) \geq m$ for all $l.$
Then define
\[
p_r
 = \sum_{l = 1}^t \sum_{j = (r - 1) d (l) + 1}^{r d (l)} e_{j, j}^{(l)}
\]
for $1 \leq r \leq m,$ and
\[
p_0
 = \sum_{l = 1}^t \sum_{j = m d (l) + 1}^{z (l)} e_{j, j}^{(l)}.
\]
The commutator estimates follow because $p_r \in E$
for $0 \leq r \leq m,$
and all the remaining statements are clear.
\end{proof}

\begin{lem}\label{LemForRPImpTRP}
Let $A$ be a stably finite simple unital \ca\   %
and let $\af \in \Aut (A).$
Assume either that $A$ has tracial rank zero,
or that $A$ is approximately divisible,
every quasitrace on $A$ is a trace,
and \pj s in $A$ distinguish the tracial states of $A.$
Suppose that there is an unbounded subset $S \subset \N$
such that for every finite set $F \subset A,$ every $\ep > 0,$
and every $n \in S,$
there are \mops\  %
\[
e_0, \, e_1, \, \ldots, \, e_{n - 1},
       \, f_0, \, f_1, \, \ldots, \, f_n \in A
\]
satisfying:
\begin{itemize}
\item[(1)]
$\| \af (e_j) - e_{j + 1} \| < \ep$ for $0 \leq j \leq n - 2$
and $\| \af (f_j) - f_{j + 1} \| < \ep$ for $0 \leq j \leq n - 1.$
\item[(2)]
$\| e_j a - a e_j \| < \ep$ for $0 \leq j \leq n - 1$ and all $a \in F,$
and
$\| f_j a - a f_j \| < \ep$ for $0 \leq j \leq n$ and all $a \in F.$
\item[(3)]
$\ta \left( 1 - \sum_{j = 0}^{n - 1} e_j + \sum_{j = 0}^{n} f_j \right)
  < \ep$
for every $\ta \in T (A).$
\end{itemize}
Then $\af$ has the \tRp.
\end{lem}

\begin{proof}
If $A$ has tracial rank zero,
then ${\mathrm{RR}} (A) = 0$ by Theorem~3.4 of~\cite{LnTAF},
and the order on projections over $A$ is determined by traces
by Theorems~5.8 and~6.8 of~\cite{LnTTR}.
Under the other hypotheses, these conclusions follow
from Corollary~3.9(b) and Theorem~1.4(e) of~\cite{BKR}.
Accordingly, we verify the conditions of Lemma~\ref{DfnUsingTrace}.

Let $F \subset A$ be a finite set, let $\ep > 0,$ and let $n \in \N.$
\Wolog\  $\| a \| \leq 1$ for every $a \in F.$
Choose $m \in \N,$ of the form $m = 2^{m_0},$
and so large that $\frac{1}{m} < \frac{1}{3} \ep.$
Choose $N \in S$ with $N > n [m (n + 1) + 1].$
Set
\[
\ep_0 = \frac{\ep}{4 (N + 2) m}.
\]
Choose $\ep_1 > 0$ with
\[
\ep_1 \leq \min \left( \frac{\ep_0}{2}, \, \frac{\ep}{3},
        \, \frac{\ep}{2 (2 N + 2)^2} \right),
\]
and also so small that whenever $B$ is a unital \ca,
whenever
\[
e_1, \, e_2, \, \ldots, \, e_{2 N - 1}
\andeqn
f_1, \, f_2, \, \ldots, \, f_{2 N - 1}
\]
are sets of orthogonal \pj s in $B$ with $\| e_k - f_k \| < \ep_1$
for $1 \leq k \leq 2 N - 1,$
then there is a unitary $u \in B$ such that
$\| u - 1 \| < \ep_0$ and $u e_k u^* = f_k$ for $1 \leq k \leq 2 N - 1.$

Apply the hypotheses with $N$ in place of $n,$
with $F$ as given, and with $\ep_1$ in place of $\ep.$
Let
$p_0, \, p_1, \, \ldots, \, p_{N - 1},
       \, q_0, \, q_1, \, \ldots, \, q_N \in A$
be the resulting \pj s.
Set
\[
r = 1 - \sum_{k = 0}^{N - 1} p_k - \sum_{k = 0}^{N} q_k.
\]
By the choice of $\ep_1,$ there is a unitary $u \in A$
such that $\| u - 1 \| < \ep_0,$
such that $u \af (p_k) u^* = p_{k + 1}$ for $1 \leq k \leq N - 2,$
and such that $u \af (q_k) u^* = q_{k + 1}$ for $1 \leq k \leq N - 1.$
Set $\bt = \Ad (u) \circ \af,$
giving $\bt (p_k) = p_{k + 1}$ and $\bt (q_k) = q_{k + 1}$
for appropriate $k.$

For $a \in A$ define
\[
E (a) = r a r + \sum_{k = 0}^{N - 1} p_k a p_k
                      + \sum_{k = 0}^{N} q_k a q_k.
\]
If $a \in F,$ then we can write $a$ as a sum of $(2 N + 2)^2$ terms
of the form $r a p_k,$ $p_j a p_k,$ etc.,
of which $2 N + 2$ appear in the formula for $E (a)$
and all the rest have norm dominated by
$\max_k \| [ p_k, a] \|$ or $\max_k \| [ q_k, a] \|.$
Accordingly, $\| E (a) - a \| < (2 N + 2)^2 \ep_1 \leq \frac{1}{2} \ep.$

We now carry out a construction involving the \pj s
$p_0, \, p_1, \, \ldots, \, p_{N - 1}.$
We do the same thing with $q_0, \, q_1, \, \ldots, \, q_N,$
but only describe the outcome afterwards.

Set
\[
F_0 = \bigcup_{k = 0}^{N - 1} \{ \bt^{-k} (p_k a p_k) \colon a \in F \}.
\]
Use Lemma~\ref{ZhangTAF} or Lemma~\ref{ZhangAppDiv},
depending on what we are assuming about $A,$
to write $p_0$ as a sum of orthogonal \pj s,
\[
p_0 = p_{0, 0} + p_{0, 1} + \cdots + p_{0, m}
\]
with
\[
p_{0, 0} \precsim p_{0, 1} \sim p_{0, 2} \sim \cdots \sim p_{0, m},
\]
and such that $\| [ p_{0, j}, b] \| < \ep_0$ for $0 \leq j \leq m$
and $b \in F_0.$
For $1 \leq k \leq N - 1$ and $0 \leq j \leq m,$
set $p_{k, j} = \bt^k (p_{0, j}) \leq p_k.$

We require estimates involving the $p_{k, j}.$
First,
\[
\| \af (p_{k, j}) - p_{k + 1, \, j} \| \leq 2 \| u - 1 \| < 2 \ep_0.
\]
Second,
we claim that if $a \in F$ then $\| [p_{k, j}, a] \| < 2 \ep_0.$
To see this, write
\begin{align*}
\| [p_{k, j}, a] \|
 & \leq \| p_{k, j} \| \cdot \| p_k a - p_k a p_k \|
            + \| [ p_k a p_k, \, p_{k, j}] \|
            + \| p_k a p_k - a p_k \| \cdot \| p_{k, j} \|  \\
 & \leq \| [p_k, a] \| + \| [ \bt^{-k} (p_k a p_k), \, p_{0, j} ] \|
            + \| [p_k, a] \|  \\
 & < \ep_1 + \ep_0 + \ep_1
   \leq 2 \ep_0.
\end{align*}
This proves the claim.

Set $N_0 = m (n + 1) + 1.$
We define subsets
\[
Y, I_0, I_1, \ldots, I_n
 \subset \{ 0, 1, \ldots, N_0 - 1 \} \times \{ 0, 1, \ldots, m \},
\]
which form a partition of this set, as follows.
Set
\[
Y = \{ (0, 0), \, (n + 1, \, 1), \, \ldots, \, (m (n + 1), \, m) \}.
\]
For $0 \leq j \leq m$ define
\begin{align*}
I_0^{(j)}
& = \{ (0, j), \, (n + 1, \, j), \, \ldots, \, ((j - 1) (n + 1), \, j),
               \\
& \hspace*{4em}    \, (j (n + 1) + 1, \, j),
     \, \ldots, \, ((m - 1) (n + 1) + 1, \, j) \}.
\end{align*}
Thus,
\[
I_0^{(0)}
  = \{ (1, 0), \, (n + 2, \, 0),
       \, \ldots, \, ((m - 1) (n + 1) + 1, \, 0) \}
\]
and
\[
I_0^{(m)}
  = \{ (0, m), \, (n + 1, \, m),
       \, \ldots, \, ((m - 1) (n + 1), \, m) \}.
\]
Then set
\[
I_0 = I_0^{(0)} \cup I_0^{(1)} \cup \cdots \cup I_0^{(m)}
\]
and for $1 \leq l \leq n$ set
\[
I_l = \{ (k + l, \, j ) \colon (k, j) \in I_0 \}.
\]
There is one more important property:
for $0 \leq k \leq N_0 - 1,$
there is at most one $j$ such that $(k, j) \in Y.$

Now write $N = d (n + 1) + s$ with $0 \leq s \leq n.$
Set $d_0 = d - s m.$
The condition on $N$ guarantees that $d_0 \geq 0.$
We define subsets
\[
Z, L_0, L_1, \ldots, L_n
 \subset \{ 0, 1, \ldots, N - 1 \} \times \{ 0, 1, \ldots, m \},
\]
which form a partition of this set, as follows.
Set
\[
Z = \{ (k + t N_0, \, j)
 \colon {\mbox{$(k, j) \in Y$ and $0 \leq t \leq s - 1$}} \}.
\]
Set
\[
J_l = \{ (k + t N_0, \, j)
 \colon {\mbox{$(k, j) \in I_l$ and $0 \leq t \leq s - 1$}} \}.
\]
These sets form a partition of
\[
\{ 0, 1, \ldots, s N_0 - 1 \} \times \{ 0, 1, \ldots, m \}.
\]
Further set
\[
K_0 = \{ (s N_0 + t (n + 1), \, j)
   \colon {\mbox{$0 \leq t \leq d_0 - 1$ and $0 \leq j \leq m$}} \}
\]
and
\[
K_l = \{ (k + l, \, j ) \colon (k, j) \in K_0 \}
\]
for $1 \leq l \leq n.$
Then set $L_l = J_l \cup K_l.$
Note that $L_l = \{ (k + l, \, j ) \colon (k, j) \in L_0 \}.$

We now introduce the notation $p_T = \sum_{(k, j) \in T} p_{k, j}$
for any subset
\[
T \subset \{ 0, 1, \ldots, N - 1 \} \times \{ 0, 1, \ldots, m \}.
\]
Define
\[
f_0 = p_{L_0}, \,\,\,\,\,\,
f_1 = p_{L_1}, \,\,\,\,\,\,
\ldots, \,\,\,\,\,\,
f_n = p_{L_n}, \andeqn
f = p_Z.
\]
These are orthogonal \pj s which add up to $\sum_{k = 0}^{N - 1} p_k.$
For $0 \leq l \leq n - 1,$ we have
\[
\| \af (f_l) - f_{l + 1} \|
  \leq \sum_{(k, j) \in L_0}
         \| \af (p_{k + l, \, j}) - p_{k + l + 1, \, j} \|
  < 2 \card (L_0) \ep_0
  \leq 2 N m \ep_0,
\]
and for $0 \leq l \leq n$ and $a \in F$ we have
\[
\| [f_l, a] \|
  \leq \sum_{(k, j) \in L_0} \| [p_{k + l, \, j}, \, a] \|
  < 2 \card (L_0) \ep_0
  \leq 2 N m \ep_0.
\]
Furthermore, for $\ta \in T (A)$ we can estimate $\ta (f)$
as follows.
For $0 \leq k \leq N - 1,$ there is at most one $j$ with $(k, j) \in Z.$
We have
\[
p_{k, 0} \precsim p_{k, 1} \sim p_{k, 2} \sim \cdots \sim p_{k, m},
\]
so that $\ta (p_{k, j}) \leq \frac{1}{m} \ta (p_k).$
Therefore
\[
\ta (f) \leq \frac{1}{m} \sum_{k = 0}^{N - 1} \ta (p_k)
  \leq \frac{1}{m} < \frac{\ep}{3}.
\]

Applying the same construction to $q_0, \, q_1, \, \ldots, \, q_N,$
we obtain orthogonal \pj s
$g, g_0, g_1, \ldots, g_n$ which add up to $\sum_{k = 0}^{N} q_k,$
and such that
$\| \af (g_l) - g_{l + 1} \| < 2 (N + 1) m \ep_0$
for $0 \leq l \leq n - 1,$
such that $\| [g_l, a] \| < 2 (N + 1) m \ep_0$
for $0 \leq l \leq n$ and $a \in F,$
and such that $\ta (g) < \frac{1}{3} \ep$
for $\ta \in T (A).$

Now set $e_l = f_l + g_l$ for $0 \leq l \leq n,$
and set $e = f + g + r = 1 - \sum_{l = 0}^n e_l.$
This gives
\[
\| \af (e_l) - e_{l + 1} \| < 2 N m \ep_0 + 2 (N + 1) m \ep_0 \leq \ep
\]
for $0 \leq l \leq n - 1,$
\[
\| [e_l, a] \| < 2 N m \ep_0 + 2 (N + 1) m \ep_0 \leq \ep
\]
for $0 \leq l \leq n$ and $a \in F,$
and $\ta (e) = \ta (f) + \ta (g) + \ta (r) < \ep$
for $\ta \in T (A).$
\end{proof}

As corollaries, we obtain the next two results.
The main difference between the first and Lemma~\ref{FOrdOnTr}
is that we do not assume that
$\ta \mapsto \ta \circ \af$ has finite order,
but we require more of the algebra.

\begin{prp}\label{UnbddSetImpTRP}
Let $A$ be a stably finite simple unital \ca\   %
and let $\af \in \Aut (A).$
Assume either that $A$ has tracial rank zero,
or that $A$ is approximately divisible,
every quasitrace on $A$ is a trace,
and that \pj s in $A$ distinguish the tracial states of $A.$
Suppose that there is an unbounded subset $S \subset \N$
such that for every finite set $F \subset A,$ every $\ep > 0,$
and every $n \in S,$
there are \mops\  %
$e_0, \, e_1, \, \ldots, \, e_n \in A$
satisfying:
\begin{itemize}
\item[(1)]
$\| \af (e_j) - e_{j + 1} \| < \ep$ for $0 \leq j \leq n - 1.$
\item[(2)]
$\| e_j a - a e_j \| < \ep$ for $0 \leq j \leq n$ and all $a \in F.$
\item[(3)]
$\ta \left( 1 - \sum_{j = 0}^{n} e_j \right) < \ep$
for every $\ta \in T (A).$
\end{itemize}
Then $\af$ has the \tRp.
\end{prp}

\begin{proof}
This is the special case of Lemma~\ref{LemForRPImpTRP}
in which always $e_j = 0$ for all $j.$
\end{proof}

\begin{thm}\label{RPImpTRP}
Let $A$ be a stably finite simple unital \ca\   %
and let $\af \in \Aut (A).$
Assume either that $A$ has tracial rank zero,
or that $A$ is approximately divisible,
every quasitrace on $A$ is a trace,
and that \pj s in $A$ distinguish the tracial states of $A.$
Suppose that $\af$ has the Rokhlin property in the sense of
Definition~\ref{StrictRP}.
Then $\af$ has the \tRp.
\end{thm}

\begin{proof}
This is the special case of Lemma~\ref{LemForRPImpTRP}
in which always
$\sum_{j = 0}^{n - 1} e_j + \sum_{j = 0}^{n} f_j = 1.$
\end{proof}

We finish this section by giving several elementary consequences of the
\tRp.

\begin{lem}\label{TRPImpOuter}
Let $A$ be a \sfsuca\   and let $\af \in \Aut (A)$ have the \tRp.
Then $\af^n$ is outer for all $n \neq 0.$
\end{lem}

\begin{proof}
It suffices to consider $n > 0.$
Let $n > 0$ and let $u \in A$ be unitary;
we show $\af^n \neq \Ad (u).$
We may clearly assume $A \not\cong \C.$
Apply Definition~\ref{TRPDfn} with this value of $n,$
with $\ep = \frac{1}{n + 2},$ with $F = \{ u \},$
and with some noninvertible $x.$
Let $e_0, e_1, \ldots, e_n$ be the resulting \pj s.

We claim that $e_j \neq 0$ for all $j.$
If $e_j = 0$ for some $j,$
relation~(1) in Definition~\ref{TRPDfn}
implies that $e_0 = e_1 = \cdots = e_n = 0.$
Then relation~(3) in Definition~\ref{TRPDfn} shows that
$1 = 1 - \sum_{j = 0}^n e_j$
is \mvnt\  to a \pj\  in ${\overline{x A x}}.$
Since ${\overline{x A x}}$ is a proper \hsa,
this contradicts stable finiteness,
and the claim follows.

Orthogonality now implies $\| e_n - e_0 \| = 1.$
Furthermore, we get
\[
\| \af^n (e_0) - e_n \| < n \ep = \frac{n}{n + 2},
\]
so
\[
\| \af^n (e_0) - e_0 \| > \frac{2}{n + 2}.
\]
However, by construction we have $\| e_0 u - u e_0 \| < \ep,$
so
\[
\| u e_0 u^* - e_0 \| < \ep = \frac{1}{n + 2}.
\]
Therefore $\af^n (e_0) \neq u e_0 u^*.$
\end{proof}

\begin{cor}\label{TRPImpSimple}
Let $A$ be a \sfsuca\   and let $\af \in \Aut (A)$ have the \tRp.
Then $C^* (\Z, A, \af)$ is simple.
\end{cor}

\begin{proof}
Using Lemma~\ref{TRPImpOuter}, this is immediate from
Theorem~3.1 of~\cite{Ks0}.
\end{proof}

\section{Rokhlin towers and subalgebras}\label{Sec:Tower}

\indent
In this section, we prove the basic approximation lemma
for actions with the \tRp.
Our first step is Proposition~\ref{OrdPjCrPrd}:
if $A$ has real rank zero and if
the order on \pj s over $A$ is determined by traces,
then for any crossed product $C^* (\Z, A, \af),$
the order with respect to
$C^* (\Z, A, \af)$ on \pj s over $A$ is determined by
traces on $C^* (\Z, A, \af),$
equivalently, by $\af$-invariant traces on $A.$
Since the proof works just as easily
for arbitrary countable amenable groups,
and since we intend to study actions of more general groups in
future work,
we give it in that generality.

\begin{ntn}\label{AffNtn}
For any compact convex set $\Dt$ in a topological vector space,
we let $\Aff (\Dt)$ be the set of all real valued \ct\  affine
functions on $\Dt.$
\end{ntn}

We are, of course, particularly interested in $\Aff (T (A)).$

The proof of Proposition~\ref{OrdPjCrPrd} requires two lemmas.

\begin{lem}\label{FAAP}
Let $A$ be a unital \ca,
and let $\af \colon \Gm \to \Aut (A)$ be an action of a countable
amenable group.
Let $f_1, \ldots, f_l \in \Aff (T (A))$ have the property that
$f_j (\ta) > 0$ for all $\Gm$-invariant $\ta \in T (A).$
Then there exist $n$ and $\gm_1, \ldots, \gm_n \in \Gm$
such that for all $\ta \in T (A)$ we have
\[
\frac{1}{n} \sum_{k = 1}^n f_j (\ta \circ \af_{\gm_k}^{-1} ) > 0
\]
for $1 \leq j \leq l.$
\end{lem}

\begin{proof}
The action of $\Gm$ on $T (A)$ will be denoted by
$(\gm \ta) (a) = \ta (\af_{\gm}^{-1} (a)).$

Since $\Gm$ is amenable, there exists a F{\o}lner sequence in $\Gm,$
that is, a sequence of nonempty finite subsets $F_n \subset \Gm$ such
that
\[
\limi{n} \frac{ \card (F_n \triangle \gm F_n)}{\card (F_n)} = 0
\]
for all $\gm \in \Gm.$
Define $S_n \colon T (A) \to T (A)$ by
\[
S_n (\ta) = \frac{1}{\card (F_n)} \sum_{\gm \in F_n} \gm \ta.
\]
Define
\[
Z_n = {\overline{ \scup{k = n}{\infty} S_n (T (A))}}.
\]
Then each $Z_n$ is a compact subset of $T (A),$
and $Z_1 \supset Z_2 \supset \cdots.$

We claim that if $\ta \in \bigcap_{n = 1}^{\infty} Z_n,$
then $\gm \ta = \ta$ for all $\gm \in \Gm.$
So let $\ta \in \bigcap_{n = 1}^{\infty} Z_n,$ let $\gm \in \Gm,$
let $a \in A,$ and let $\ep > 0.$
Choose $N$ so large that if $n \geq N$ then
\[
\frac{ \card (F_n \triangle \gm F_n)}{\card (F_n)}
   < \frac{\ep}{3 \| a \|}.
\]
By the definition of the weak* topology, there is
$\sm \in \bigcup_{n = N}^{\infty} S_n (T (A))$ such that
\[
| \sm (a) - \ta (a) | < \ts{ \frac{1}{3} } \ep
\andeqn
| \sm (\af_{\gm}^{-1} (a)) - \ta (\af_{\gm}^{-1} (a)) |
   < \ts{ \frac{1}{3} } \ep.
\]
Write $\sm = S_n (\rh)$ for some $n \geq N$ and $\rh \in T (A).$
Then
\begin{align*}
| \ta (\af_{\gm}^{-1} (a)) - \ta (a) |
 & < \ts{ \frac{2}{3} } \ep + | \sm (\af_{\gm}^{-1} (a)) - \sm (a) |
          \\
 & \! = \frac{2 \ep}{3}
      + \frac{1}{\card (F_n)}
         \left| \susum{\et \in F_n}
                \rh (\af_{\et}^{-1} \circ \af_{\gm}^{-1} (a))
               - \susum{\et \in F_n} \rh (\af_{\et}^{-1} (a)) \right|
          \\
 & \! \leq \frac{2 \ep}{3}
      + \frac{ \card (F_n \triangle \gm F_n) \| a \|}{\card (F_n)}
      < \ep.
\end{align*}
Since $\ep > 0$ is arbitrary, it follows that $\gm \ta = \ta.$

Now set
\[
Y_n = Z_n \cap \{ \ta \in T (A) \colon
     {\mbox{$f_j (\ta) \leq 0$ for some $j$ with $1 \leq j \leq l$}} \}.
\]
Then each $Y_n$ is compact,
and $Y_1 \supset Y_2 \supset \cdots.$
Moreover, $\bigcap_{n = 1}^{\infty} Y_n = \varnothing,$
because any element $\ta$ of this set is an invariant tracial state
such that $f_j (\ta) \leq 0$ for some $j.$
Therefore there is $n$ such that $Y_n = \varnothing.$
Now $f_j (S_n (\ta)) > 0$
for $1 \leq j \leq l$ and all $\ta \in T (A).$
Since each $f_j$ is affine, we have
\[
\frac{1}{\card (F_n)} \sum_{\gm \in F_n} f_j ( \gm \ta)
  = f_j (S_n (\ta)) > 0,
\]
as required.
\end{proof}

The following lemma is a more flexible version of a result of
Zhang~\cite{Zh7}.
In Zhang's version, which we use in the proof,
the integer $n$ of the hypotheses is required to be a power of~$2.$

\begin{lem}\label{FlexZhang}
Let $A$ be a simple unital infinite dimensional \ca\  %
with real rank zero.
Let $p \in A$ be a \pj, and let $n \in \N.$
Then there exist \pj s $p_0, p_1, \ldots, p_n \in A$ such that
\[
\sum_{k = 0}^n p_k = p, \,\,\,\,\,\,
p_1 \sim p_2 \sim \cdots \sim p_n,
\andeqn
p_0 \precsim p_1.
\]
\end{lem}

\begin{proof}
Choose $m \in \N$ such that $2^m > n^2.$
Set $N = 2^m,$ and write $N = l n + r$ for integers $r$ and $l$
such that $0 \leq r < n.$
Note that $l \geq n.$
Apply Theorem~1.1(a) of~\cite{Zh7}, obtaining
\pj s $e_0, e_1, \ldots, e_N \in A$ such that
\[
\sum_{k = 0}^N e_k = p, \,\,\,\,\,\,
e_1 \sim e_2 \sim \cdots \sim e_N,
\andeqn
e_0 \precsim e_1.
\]
Define
\[
p_0 = e_0 + e_{n l + 1} + \cdots + e_N
\]
and, for $1 \leq k \leq n,$ define
\[
p_k = e_{(k - 1) l + 1} + e_{(k - 1) l + 2} + \cdots + e_{k l}.
\]
The conditions $\sum_{k = 0}^n p_k = p$ and
$p_1 \sim p_2 \sim \cdots \sim p_n$ in the conclusion are obvious,
and $p_0 \precsim p_1$ follows from $e_0 \precsim e_{(k - 1) l + 1}$
and the fact that there are $r + 1 \leq n \leq l$ terms in the
sum defining $p_0.$
\end{proof}

\begin{prp}\label{OrdPjCrPrd}
Let $A$ be a simple unital infinite dimensional \ca\  %
with real rank zero,
and assume that the order on \pj s over $A$ is determined by traces.
Let $\af \colon \Gm \to \Aut (A)$ be an action of a countable
amenable group.
Let $p, \, q \in \Mi (A)$ be \pj s such that
$\ta (p) < \ta (q)$ for every $\Gm$-invariant tracial
state $\ta$ on $A.$
(We extend $\ta$ to $\Mi (A)$ in the obvious way.)
Then there is $s \in \Mi ( C^* (\Gm, A, \af))$
such that
\[
s s^* = p, \,\,\,\,\,\, s s^* \leq q, \andeqn s s^* \in \Mi (A).
\]
In particular, $p \precsim q$ in $\Mi ( C^* (\Gm, A, \af)).$
\end{prp}

\begin{proof}
Throughout the proof, we regard elements of $T (A)$ as being
defined on all of $\Mi (A)$ in the obvious way.

Define $f \in \Aff (T (A))$ by $f (\ta) = \ta (q) - \ta (p).$
Use Lemma~\ref{FAAP} on this function $f$ to find
$n \in \N$ and $\gm_1, \ldots, \gm_n \in \Gm$
such that for all $\ta \in T (A)$ we have
\[
g (\ta) = \frac{1}{n} \sum_{k = 1}^n [\ta ( \af_{\gm_k}^{-1} (q) )
                    - \ta ( \af_{\gm_k}^{-1} (p) )] > 0.
\]
Then set $\ep = \inf_{\ta \in T (A)} g (\ta),$
which is strictly positive because $T (A)$ is compact.
Also set $M = \sup_{\ta \in T (A)} \ta (q),$
which is finite for the same reason.
Choose $N \in \N$ such that
\[
\frac{M}{N n} < \frac{\ep}{3}.
\]
Use Lemma~\ref{FlexZhang} on $p$ with $N n - 1$ in place of $n,$
calling the resulting \pj s $p_0, p_1, \ldots, p_{N n - 1},$
and on $q$ with $N n$ in place of $n,$
calling the resulting \pj s $q_0, q_1, \ldots, q_{N n}.$

We now claim that
\[
\sum_{k = 1}^n \ta ( \af_{\gm_k}^{-1} (p_1) )
 < \sum_{k = 1}^n \ta ( \af_{\gm_k}^{-1} (q_1) )
\]
for all $\ta \in T (A).$
To see this, use
\[
p_1 \sim p_2 \sim \cdots \sim p_{N n - 1}
\andeqn
q_0 \precsim q_1 \sim q_2 \sim \cdots \sim q_{N n}
\]
to get
\[
\frac{1}{n} \sum_{k = 1}^n \ta ( \af_{\gm_k}^{-1} (p) )
  \geq \frac{N n - 1}{n} \sum_{k = 1}^n \ta ( \af_{\gm_k}^{-1} (p_1) )
\]
and
\[
\frac{1}{n} \sum_{k = 1}^n \ta ( \af_{\gm_k}^{-1} (q) )
  \leq \frac{N n + 1}{n} \sum_{k = 1}^n \ta ( \af_{\gm_k}^{-1} (q_1) ).
\]
We also have
\[
\frac{2}{n} \sum_{k = 1}^n \ta ( \af_{\gm_k}^{-1} (q_1) )
  \leq \frac{2}{n} \sum_{k = 1}^n
            \frac{\ta ( \af_{\gm_k}^{-1} (q) )}{N n}
  \leq \frac{2 M}{N n} < \frac{2 \ep}{3}.
\]
Using this result at the last step, we get
\begin{align*}
\frac{N n - 1}{n} \sum_{k = 1}^n \ta ( \af_{\gm_k}^{-1} (p_1) )
 & \leq \frac{1}{n} \sum_{k = 1}^n \ta ( \af_{\gm_k}^{-1} (p) )
   \leq \frac{1}{n} \sum_{k = 1}^n \ta ( \af_{\gm_k}^{-1} (q) )
                 - \ep  \\
 & \hspace*{-1em}
   \leq \frac{N n + 1}{n} \sum_{k = 1}^n \ta ( \af_{\gm_k}^{-1} (q_1) )
                 - \ep
   \leq \frac{N n - 1}{n} \sum_{k = 1}^n \ta ( \af_{\gm_k}^{-1} (q_1) )
                 - \frac{\ep}{3}.
\end{align*}
The claim follows by dividing by $(N n - 1) / n.$

Now for $1 \leq k \leq n$ define
\[
e_k = p_{(k - 1) N} + p_{(k - 1) N + 1} + \cdots + p_{k N - 1}
\]
and
\[
f_k = q_{(k - 1) N + 1} + q_{(k - 1) N + 2} + \cdots + q_{k N}.
\]
(We do not use $q_0.$)
Regarding $e$ and $f$ as elements of $M_r (A)$ for suitable $r,$
further define
\[
e = \diag (e_1, \ldots, e_n)
\andeqn {\overline{e}}
  = \diag (\af_{\gm_1} (e_1), \ldots, \af_{\gm_n} (e_n) ),
\]
and
\[
f = \diag (f_1, \ldots, f_n)
\andeqn {\overline{f}}
  = \diag (\af_{\gm_1} (f_1), \ldots, \af_{\gm_n} (f_n) ),
\]
which are all \pj s in $M_{r n} (A).$
By construction, in $\Mi (A)$ we have
\[
p \sim e \andeqn f \sim \sum_{m = 1}^{N n} q_m = q - q_0.
\]

For $\gm \in \Gm$
let $u_{\gm}$ be the standard unitary in $C^* (\Gm, A, \af)$
which implements $\af_{\gm}.$
Set
\[
v = \diag (1_r \otimes u_{\gm_1},
             \, \, \ldots, \, 1_r \otimes u_{\gm_n})
 \in M_{r n} (A),
\]
so that $v e v^* = {\overline{e}}$ and $v f v^* = {\overline{f}}.$
The claim proved above implies that
$\ta ( {\overline{e}} ) < \ta ( {\overline{f}} )$
for all $\ta \in T (A),$
whence ${\overline{e}} \precsim {\overline{f}}$ in $M_{r n} (A).$

We now have enough to get
$p \precsim q$ in $M_{\infty} (C^* (\Gm, A, \af)),$
but we need more to get the stronger statement in the conclusion.
Since $M_{r n} (A)$ has real rank zero,
Theorem~1.1 of~\cite{Zh} implies that \pj s in $M_{r n} (A)$ satisfy
Riesz decomposition,
so there are \pj s $g_1, \ldots, g_n \in M_{r} (A)$ such that
$g_k \leq \af_{\gm_k} (f_k)$ for all $k$ and,
with $g = \diag (g_1, \ldots, g_n),$
we have ${\overline{e}} \sim g$ in $M_{r n} (A).$
Then
\[
{\overline{g}} = v^* g v
 = \diag (\af_{\gm_1}^{-1} (g_1), \ldots, \af_{\gm_n}^{-1} (g_n) )
 \in M_{r n} (A)
\]
satisfies ${\overline{g}} \sim g$ in $C^* (\Gm, A, \af)$
and ${\overline{g}} \leq f.$
Since $f \sim q - q_0$ in $\Mi (A),$
there is a \pj\  $h \in M_r (A)$
such that $h \leq q - q_0$ and ${\overline{g}} \sim h$ in $\Mi (A).$
Thus, in $M_{\infty} (C^* (\Gm, A, \af))$ we have
\[
p \sim e \sim {\overline{e}} \sim g
   \sim {\overline{g}} \sim h \leq q - q_0 \leq q,
\]
with $h \in M_r (A).$
\end{proof}

\begin{lem}\label{Basic1}
Let $A$ be a \sfsuca\  with real rank zero such that the order on
\pj s over $A$ is determined by traces.
Let $\af \in \Aut (A)$ have the \tRp.
Let $\io \colon A \to C^* (\Z, A, \af)$ being the inclusion map.
Then for every finite set $F \subset C^* (\Z, A, \af),$
every $\ep > 0,$
every $N \in \N,$
every nonzero positive element $z \in C^* (\Z, A, \af),$
and every sufficiently large $n \in \N$
(depending on $F,$ $\ep,$ $N,$ and $z$),
there exist a \pj\  $e \in A \subset C^* (\Z, A, \af),$
a unital subalgebra $D \subset e C^* (\Z, A, \af) e,$
a \pj\  $p \in D,$
a \pj\  $f \in A,$ and an isomorphism
$\ph \colon M_n \otimes f A f \to D,$ such that:
\begin{itemize}
\item[(1)]
With $( e_{j, k} )$ being the standard system of matrix units for $M_n,$
we have $\ph (e_{1, 1} \otimes a) = \io (a)$ for all $a \in f A f$
and $\ph (e_{k, k} \otimes 1) \in \io (A)$ for $1 \leq k \leq n.$
\item[(2)]
With $( e_{j, k} )$ as in~(1), we have
$\| \ph (e_{j, j} \otimes a) - \af^{j - 1} ( \io (a)) \|
           \leq \ep \| a \|$
for all $a \in f A f.$
\item[(3)]
For every $a \in F$ there exist $b_1, \, b_2 \in D$ such that
$\| p a - b_1 \| < \ep,$ $\| a p - b_2 \| < \ep,$
and $\| b_1 \|, \, \| b_2 \| \leq \| a \|.$
\item[(4)]
There is $m \in \N$ such that $2 m / n < \ep$ and
$p = \sum_{j = m + 1}^{n - m} \ph (e_{j, j} \otimes 1).$
\item[(5)]
The \pj\  $1 - p$ is \mvnt\  in $C^* (\Z, A, \af)$ to a
\pj\  in the \hsa\  of $C^* (\Z, A, \af)$ generated by $z.$
\item[(6)]
There are $N$
\mops\  $f_1, f_2, \ldots, f_N \in p D p,$
each of which is \mvnt\  in $C^* (\Z, A, \af)$ to $1 - p.$
\end{itemize}
\end{lem}

\begin{proof}
We first make a simplification:
We need not check the estimates
$\| b_1 \|, \, \| b_2 \| \leq \| a \|$
in Condition~(3) of the conclusion.
To prove this, \wolog\  $\| a \| \leq 1$ for all $a \in F.$
Apply the weaker statement with $\frac{1}{2} \ep$
in place of $\ep,$ and with all other parameters the same.
Let $c_1$ and $c_2$ be the resulting elements
in Condition~(3) of the conclusion.
Then $\| c_1 \|, \, \| c_2 \| \leq 1 + \ts{\frac{1}{2}} \ep.$
Set
\[
b_1 = \left( \frac{1}{1 + \ts{ \frac{1}{2} } \ep} \right) c_1 \andeqn
b_2 = \left( \frac{1}{1 + \ts{ \frac{1}{2} } \ep} \right) c_2.
\]
One checks that $\| b_1 - c_1 \|\leq \frac{1}{2} \ep,$
so $\| b_1 - p a \| < \ep.$
Similarly $\| b_2 - a p \| < \ep.$
This proves the reduction.

Now we do the main part of the proof.
Let $\ep > 0,$ and let $F \subset C^* (\Z, A, \af)$ be a finite set.
Let $N \in \N,$ and let $z \in C^* (\Z, A, \af)$ be a nonzero
positive element.

Let $u$ be the standard unitary
in the crossed product $C^* (\Z, A, \af).$
We regard $A$ as a subalgebra of $C^* (\Z, A, \af)$ in the usual way.
Choose $m \in \N$ such that for every $x \in F$ there are
$a_l \in A$ for $-m \leq l \leq m$ such that
\[
\left\| x - \ssum{l = -m}{m} a_l u^l \right\| < \frac{\ep}{2}.
\]
For each $x \in F$ choose one such expression,
and let $S \subset A$ be a finite set which contains all the
coefficients used for all elements of $F.$
Let $M = 1 + \sup_{a \in S} \| a \|.$

Since $A$ has Property~(SP),
and since (by Lemma~\ref{TRPImpOuter})
all nontrivial powers of $\af$ are outer,
we can apply Theorem~4.2 of \cite{JO}, with $N = \{ 1 \},$
to find a nonzero \pj\  $q \in A$ which is \mvnt\  in
$C^* (\Z, A, \af)$ to a \pj\  in
${\overline{z C^* (\Z, A, \af) z}}.$
Moreover, Lemma~\ref{FlexZhang} provides
nonzero orthogonal \mvnt\   \pj s
$g_0, \, g_1, \, \ldots, \, g_{2 m} \in q A q.$

Since $A$ is simple, $g_0$ is a nonzero \pj, and the tracial state
space $T (A)$ of $A$ is weak* compact, we have
$\dt = \inf_{\ta \in T (A)} \ta (g_0) > 0.$
Now let $n \in \N$ be any integer such that
\[
n > \max \left(
     \frac{1}{\dt}, \, (N + 2) (2 m + 1), \, \frac{2 m}{\ep} \right).
\]
Set
\[
\ep_0 = \frac{\ep}{10 (2 m + 1) n^2 M}.
\]

Choose $\ep_1 > 0$ so small that whenever
$e_1, e_2, \ldots, e_n$ are \mops\  in a unital \ca\  $B$
and $u \in B$ is a unitary such that
$\| u e_j u^* - e_{j + 1} \| < \ep_1$
for $1 \leq j \leq n,$
then there is a unitary $v \in B$ such that $\| v - u \| < \ep_0$
and $v e_j v^* = e_{j + 1}$ for $1 \leq j \leq n.$
Further use Lemma~\ref{FlexZhang} to find nonzero orthogonal \mvnt\  %
\pj s $h_1, \, h_2, \, \ldots, \, h_{n + 2} \leq g_0.$

Apply the \tRp\  (Definition~\ref{TRPDfn}) with
$n - 1$ in place of $n,$
with $\min (1, \ep_1, \ep_0)$ in place of $\ep,$
with $S$ in place of $F,$ and with $h_1$ in place of $x.$
Call the resulting \pj s $e_1, e_2, \ldots, e_n,$
and let $e = \sum_{j = 1}^{n} e_j.$
Apply the choice of $\ep_1$ to these \pj s and the
standard unitary $u,$ obtaining a unitary $v \in C^* (\Z, A, \af)$
as in the previous paragraph.

Set $f = e_1.$
The elements $e_j v^{j - k} e_k,$
for $1 \leq j, \, k \leq n,$
can be seen to satisfy the relations for matrix units $e_{j, k}.$
So there is a unique injective \hm\  %
$\ph \colon M_{n} \otimes f A f \to C^* (\Z, A, \af)$
such that $\ph (e_{1, 1} \otimes a) = a$ for $a \in f A f$
and $\ph (e_{j, k} \otimes f) = e_j v^{j - k} e_k$
for $1 \leq j, \, k \leq n.$
Let $D$ be the range of $\ph,$
so that $\ph \colon M_{n} \otimes f A f \to D$ is an isomorphism,
as required.
Condition~(1) of the conclusion is immediate.
For $1 \leq j, \, k \leq n$ and $a \in f A f$ we have
\[
\ph (e_{j, k} \otimes a)
  = \ph (e_{j, 1} \otimes 1) \ph (e_{1, 1} \otimes a)
          \ph (e_{1, k} \otimes 1)
  = e_j v^{j - 1} a v^{1 - k} e_k = v^{j - 1} a v^{1 - k}.
\]
In particular, if $j = k$ then
\begin{align*}
\| \ph (e_{j, j} \otimes a) - \af^{j - 1} (a) \|
&  \leq 2\| a \| \cdot \| v^{j - 1} - u^{j - 1} \|
   \leq 2\| a \| \cdot (j - 1) \| v - u \|       \\
&  \leq 2n \ep_0 \| a \| \leq \ep \| a \|.
\end{align*}
This is Condition~(2) of the conclusion.

Let $p = \sum_{j = m + 1}^{n - m} e_j,$
and note that
$\sum_{j = m + 1}^{n - m} \ph (e_{j, j} \otimes 1) = p.$
Condition~(4) of the conclusion now follows from the choice of $n.$

We now claim that if $y = \sum_{l = -m}^{m} a_l v^l$
with $a_l \in A$ for $-m \leq l \leq m,$
and if $[e_j, a_l] = 0$
for $-m \leq l \leq m$ and $1 \leq j \leq n,$
then there are $d_1, \, d_2 \in D$ such that
\[
\| p y - d_1 \|, \, \| y p - d_2 \| < 2 M n (n - 2 m) (2 m + 1) \ep_0.
\]
We produce $d_1$; the proof for $d_2$ is essentially the same.
We write
\[
p y = \sum_{j= m + 1}^{n - m} \sum_{l = -m}^{m} e_j a_l v^l
   = \sum_{j= m + 1}^{n - m}
            \sum_{l = -m}^{m} (e_j a_l e_j) (e_j v^l e_{j - l}).
\]
Since $v^{j - 1} e_1 v^{- j + 1} = e_j,$ we have
\begin{align*}
\| \ph (e_{j, j} \otimes f \af^{- j + 1} (a_{l} ) f )
          - e_{j} a_{l} e_{j} \|
 & = \| e_j v^{j - 1} e_1 u^{- j + 1} a_l u^{j - 1} e_1 v^{- j + 1} e_j
                    - e_{j} a_{l} e_{j} \|  \\
 & \leq 2 \| a_{l} \| \cdot \| u^{j - 1} - v^{j - 1} \|  \\
 & < 2 M (j - 1) \ep_0
   \leq 2 M n \ep_0,
\end{align*}
so
\[
\| \ph (e_{j, \, j - l} \otimes f \af^{- j + 1} (a_{l}) f )
    - (e_{j} a_{l} e_{j}) (e_{j} v^{l} e_{j - l} ) \| < 2 M n \ep_0.
\]
Therefore
\[
\sum_{j= m + 1}^{n - m} \sum_{l = -m}^{m}
         (e_{j} a_{l} e_{j}) (e_{j} v^{l} e_{j - l})
\]
differs from an element of $D$ by less than
$2 M n (n - 2 m) (2 m + 1) \ep_0.$
The claim follows.

We next prove Condition~(3) of the conclusion.
Let $x \in F.$
Choose
\[
b_{-m}, \, b_{- m + 1}, \, \ldots, b_m \in S
\]
such that
\[
\left\| x - \ssum{l = -m}{m} b_l u^l \right\| < \frac{\ep}{2}.
\]
For $-m \leq l \leq m$ define
\[
a_l = (1 - e) b_l (1 - e) + \sum_{j = 1}^{n} e_j b_l e_j.
\]
We write
\[
b_l - a_l = \sum_{j = 1}^{n} [e_j a_l (1 - e) + (1 - e) a_l e_j]
     + \sum_{i = 1}^{n} \sum_{j \neq i} e_i a_l e_j,
\]
so that the estimate
$\| [a_l, e_j] \| < \ep_1 \leq \ep_0$ implies
\[
\| b_l - a_l \| < [2 n + n (n - 1)] \ep_0 < 2 n^2 \ep_0.
\]
Moreover, from $\| v - u \| < \ep_0$
we get $\| v^l - u^l \| < m \ep_0$ for $-m \leq l \leq m.$
Therefore, with $y = \sum_{l = -m}^{m} a_l v^l,$ we get
\begin{align*}
\left\| x - y \right\|
 & \leq \left\| x - \ssum{l = -m}{m} b_l u^l \right\|
    + \sum_{l = -m}^{m} \| b_l \| \cdot \| v^l - u^l \|
    + \sum_{l = -m}^{m} \| b_l - a_l \|  \\
 & < \ts{ \frac{1}{2} } \ep + (2 m + 1) M m \ep_0
      + (2 m + 1) \cdot 2 n^2 \ep_0.
\end{align*}
According to our claim, there is $d \in D$ such that
$\| p y - d \| < 2 M n (n - 2 m) (2 m + 1) \ep_0.$
Then
\begin{align*}
\| p x - d \|
& < \ts{ \frac{1}{2} } \ep
       + (2 m + 1) [ M m + 2 n^2 + 2 M n (n - 2 m)] \ep_0  \\
& \leq \ts{ \frac{1}{2} } \ep + (2 m + 1) \cdot 5 M n^2 \ep_0
  \leq \ep.
\end{align*}
This is one half of Condition~(3) of the conclusion.
The other half is proved similarly.

It remains to verify Conditions~(5) and~(6) of the conclusion.
We have
\[
1 - p = 1 - e + \sum_{j = 1}^{m} e_j + \sum_{j = n - m + 1}^{n} e_j.
\]
By construction we have $1 - e \precsim h_1 \leq g_0.$
Now let $\ta$ be any $\af$-invariant tracial state on $A.$
Then $\ta (e_j) = \ta (e_1)$ for all $j,$
whence $\ta (e_j) \leq \frac{1}{n}.$
The inequality
\[
n > \frac{1}{\dt} \geq \frac{1}{\ta (g_0)}
\]
therefore implies $\ta (e_j) < \ta (g_0).$
Since all $g_j$ are \mvnt, it follows that for any $\af$-invariant
tracial sate $\ta$ we have
\[
\ta (e_j) < \ta (g_j) \andeqn \ta (e_{n - j}) < \ta (g_{m + j})
\]
for $1 \leq j \leq m.$
So Proposition~\ref{OrdPjCrPrd} implies that
\[
e_{j} \precsim g_j \andeqn e_{n - j} \precsim g_{m + j}
\]
in $C^* (\Z, A, \af)$ for $1 \leq j \leq m.$
Thus
\[
1 - p \precsim \sum_{j = 0}^{2 m} g_j \leq q,
\]
which is \mvnt\  in
$C^* (\Z, A, \af)$ to a \pj\  in the \hsa\  %
${\overline{z C^* (\Z, A, \af) z}}.$
This is Condition~(5) of the conclusion.

Finally, we prove~(6).
Let $\ta \in T (A)$ be $\af$-invariant.
By construction, we have $1 - e \precsim h_1.$
Since $h_1, \, h_2, \, \ldots, \, h_{n + 2} \leq g_0 \leq 1$
are orthogonal \mvnt\  \pj s,
we get $\ta (h_1) \leq (n + 2)^{-1}$ for all $\ta \in T (A),$
and in particular $\ta (1 - e) \leq (n + 2)^{-1}.$
Since all $\ta (e_j)$ are equal, we have
\[
\ta (e_j)
 \geq \frac{1}{n} \left( 1 - \frac{1}{n + 2} \right)
 > \frac{1}{n + 2}.
\]
So Proposition~\ref{OrdPjCrPrd} provides a \pj\  in
$e_j A e_j \subset e_j D e_j$ which is \mvnt\  in $C^* (\Z, A, \af)$
to $1 - e.$
Therefore, for every $k \geq 0$ with $(2 m + 1) (k + 2) \leq n,$
\[
1 - p = 1 - e + \sum_{j = 1}^{m} e_j + \sum_{j = n - m + 1}^{n} e_j
   \precsim \sum_{j = (2 m + 1) k + 1}^{(2 m + 1) (k + 1)} e_j \leq p,
\]
and the \pj\  \mvnt\  to $1 - p$ can be chosen to be in $p D p.$
Since $n \geq (N + 2) (2 m + 1),$
there are at least $N$ such \pj s.
They are orthogonal, so Condition~(6) of the conclusion is verified.
\end{proof} % of Lemma~\ref{Basic1}

Given objects satisfying part~(1) of the conclusion
of Lemma~\ref{Basic1},
we can make a useful \hm\  into $C^* (\Z, A, \af)$
which should be thought of as
a kind of twisted inclusion of $A.$

\begin{lem}\label{EmbedA}
Let $A$ be any simple unital \ca, let $\af \in \Aut (A),$
and let $\io \colon A \to C^* (\Z, A, \af)$ be the inclusion.
Let $e, f \in A$ be a \pj s, and let $n \in \N.$
Assume that there is an injective unital \hm\  %
$\ph \colon M_n \otimes f A f \to \io (e) C^* (\Z, A, \af) \io (e)$
such that,
with $( e_{j, k} )$ being the standard system of matrix units for $M_n,$
we have $\ph (e_{1, 1} \otimes a) = \io (a)$ for all $a \in f A f.$
Then there is a corner $A_0 \subset M_{n + 1} \otimes A$ which contains
\[
\left\{ \ts{ \left( \begin{array}{cc} a & 0 \\ 0 & b \end{array} \right) }
  \colon {\mbox{$a \in (1 - e) A (1 - e)$ and
                 $b \in M_n \otimes f A f$}} \right\}
\]
as a unital subalgebra, and
an injective unital \hm\  $\ps \colon A_0 \to C^* (\Z, A, \af)$
such that
\[
\ps \left( \ts{ \begin{array}{cc} a & 0 \\ 0 & b \end{array} } \right)
  = \io (a) + \ph (b)
\]
for $a \in (1 - e) A (1 - e)$ and $b \in M_n \otimes f A f.$

Moreover, for every $\af$-invariant tracial state $\ta$ on $A$
there is a tracial state $\sm$ on $C^* (\Z, A, \af)$ such that
the extension ${\overline{\ta}}$ of $\ta$ to $M_{n + 1} \otimes A$
satisfies ${\overline{\ta}} |_{A_0} = \sm \circ \ps.$
\end{lem}

\begin{proof}
Set
\[
q = \diag (1 - e, \, f, \, f, \, \ldots, \, f) \in M_{n + 1} \otimes A,
\]
and set
\[
A_0 = q (M_{n + 1} \otimes A) q
\andeqn e_0  = \diag (0, f, f, \ldots, f) \in A_0.
\]
In $M_{n + 1},$ call the matrix units $e_{j, k}$
for $0 \leq j, \, k \leq n.$
Then $q - e_0 = e_{0, 0} \otimes (1 - e).$
Define $\ps \colon A_0 \to C^* (\Z, A, \af)$ as follows.
\begin{itemize}
\item[(1)]
For $a \in (q - e_0) A_0 (q - e_0),$ write $a = e_{0, 0} \otimes x$
with $x \in (1 - e) A (1 - e),$ and set $\ps (a) = \io (x).$
\item[(2)]
For $a \in e_0 A_0 e_0,$ write
$a = \sum_{j, k = 1}^n e_{j, k} \otimes x_{j, k}$
with $x_{j, k} \in f A f$ for all $j$ and $k.$
Regard this sum as an element of $M_n \otimes f A f$ in the
obvious way, and set $\ps (a) = \ph (a).$
\item[(3)]
For $a \in (e_{j, j} \otimes f) A_0 (q - e_0)$
for some $j$ with $1 \leq j \leq n,$
write $a = e_{j, 0} \otimes x$
with $x \in f A (1 - e),$
and set $\ps (a) = \ph (e_{j, 1} \otimes f) \io (x).$
\item[(4)]
For $a \in (q - e_0) A_0 (e_{j, j} \otimes f)$
for some $j$ with $1 \leq j \leq n,$
set $\ps (a) = \ps (a^*)^*$ using~(3).
\end{itemize}
Then extend by linearity.

To prove the first part of the lemma,
it suffices to prove that $\ps$ defined this way is in fact a \hm.
It is clear that $\ps$ is linear and
that $\ps (a^*) = \ps (a)^*$ for all $a \in A_0,$
so we prove multiplicativity.
We show that $\ps (a b) = \ps (a) \ps (b)$ in four cases:
\begin{itemize}
\item[(5)]
$a \in (e_{j, j} \otimes f) A_0 (q - e_0)$
and $b \in (q - e_0) A_0 (q - e_0).$
\item[(6)]
$a \in (e_{j, j} \otimes f) A_0 (q - e_0)$
and $b \in (q - e_0) A_0 (e_{j, j} \otimes f).$
\item[(7)]
$a \in (q - e_0) A_0 (e_{j, j} \otimes f)$
and $b \in e_0 A_0 e_0.$
\item[(8)]
$a \in (q - e_0) A_0 (e_{j, j} \otimes f)$
and $b \in (e_{k, k} \otimes f) A_0 (q - e_0).$
\end{itemize}
The other $12$ cases are all of three kinds:
both $\ps (a b)$ and $\ps (a) \ps (b)$ are easily seen to be zero;
the formula $\ps (a b) = \ps (a) \ps (b)$ follows from the fact that
$\io$ is a \hm\  or $\ph$ is a \hm;
or the case follows from one of the four cases above by taking adjoints.

For~(5), write $a = e_{j, 0} \otimes x$ as in~(3)
and write $b = e_{0, 0} \otimes y$ analogously to~(1).
Then $a b = e_{j, 0} \otimes x y$ analogously to~(3), so
\[
\ps (a) \ps (b) = \ph (e_{j, 1} \otimes f) \io (x) \io (y)
 = \ph (e_{j, 1} \otimes f) \io (x y) = \ps (a b).
\]
For~(6), the analogous computation is:
$a = e_{j, 0} \otimes x,$ $b = e_{0, j} \otimes y,$ and,
using $x y \in f A f$ so that
$\io (x y) = \ph (e_{1, 1} \otimes x y),$
\begin{align*}
\ps (a) \ps (b)
& = \ph (e_{j, 1} \otimes f) \io (x) \io (y) \ph (e_{1, j} \otimes f)
  = \ph (e_{j, 1} \otimes f) \io (x y) \ph (e_{1, j} \otimes f)  \\
& = \ph (e_{j, 1} \otimes f) \ph (e_{1, 1} \otimes x y)
             \ph (e_{1, j} \otimes f)
  = \ps (a b).
\end{align*}
Similarly, in~(7) write
$a = e_{0, j} \otimes x$ with $x \in (1 - e) A f$
and $b = \sum_{j, k = 1}^n e_{j, k} \otimes y_{j, k}$
with all $y_{j, k} \in f A f$;
then
\[
a b = \sum_{k = 1}^n e_{0, k} \otimes x y_{j, k}
\]
with $x y_{j, k} \in (1 - e) A f,$
and
\begin{align*}
\ps (a) \ps (b)
& = \sum_{k = 1}^n
   \io (x) \ph (e_{1, j} \otimes f) \ph (e_{j, k} \otimes y_{j, k})
  = \sum_{k = 1}^n \io (x)
   \ph (e_{1, 1} \otimes y_{j, k}) \ph (e_{1, k} \otimes f)     \\
& = \sum_{k = 1}^n \io (x) \io (y_{j, k}) \ph (e_{1, k} \otimes f)
  = \ps (a b).
\end{align*}
Finally, in~(8) if $j \neq k$
one easily gets $\ps (a) \ps (b) = 0 = \ps (a b),$ and otherwise
one writes
$a = e_{0, j} \otimes x,$ $b = e_{j, 0} \otimes y,$ and
\begin{align*}
\ps (a) \ps (b)
& = \io (x) \ph (e_{1, j} \otimes f) \ph (e_{j, 1} \otimes f) \io (y)
  = \io (x) \ph (e_{1, 1} \otimes f) \io (y)  \\
& = \io (x) \io (f) \io (y)
  = \io (x y)
  = \ps (a b).
\end{align*}

It remains to prove the statement about the tracial states.
So let $\ta$ be an $\af$-invariant tracial state on $A.$
Let $E \colon C^* (\Z, A, \af) \to A$
be the standard conditional expectation,
and let $\sm = \ta \circ E$
be the induced tracial state on $C^* (\Z, A, \af).$
If $f = 0$ then $A_0 = A$ and $\ps = \io,$ so the statement is immediate.
Otherwise, for $a \in f A f,$ we have
\[
\sm \circ \ps (e_{1, 1} \otimes a) = \sm \circ \ph (e_{1, 1} \otimes a)
 = \sm \circ \io (a) = \ta (a).
\]
Therefore $\sm \circ \ps$ and ${\overline{\ta}}$ agree on the
full corner
$(e_{1, 1} \otimes f) (M_{n + 1} \otimes A) (e_{1, 1} \otimes f)$
of $A_0.$
So $\sm \circ \ps = {\overline{\ta}}.$
\end{proof} % of Lemma~\ref{EmbedA}

\section{Traces and order on projections in crossed
   products}\label{Sec:OrdPj}

\indent
In this section, we prove that if
$A$ is a simple unital C*-algebra with real rank zero such that
the order on projections over $A$ is
determined by traces, and if $\af \in \Aut (A)$ has the \tRp,
then the order on projections over $C^* (\Z, A, \af)$ is
determined by traces.
The methods are adapted from Section~3 of~\cite{Ph10},
and originally came from~\cite{Pt5}.
We make one small improvement.
In previous versions of this argument,
the conclusion was only that
the order on $K_0 (C^* (\Z, A, \af))$ is determined by traces,
and the result on the order on \pj s was then obtained
using stable rank one.
Here, we obtain the full result even if
$C^* (\Z, A, \af)$ does not have stable rank one.

We begin with a comparison lemma for \pj s in
crossed products by actions with the \tRp.

\begin{lem}\label{BasicOrder}
Assume the hypotheses of Lemma~\ref{EmbedA},
and assume in addition that $A$ has real rank zero
and that the order on projections over $A$ is
determined by traces.
Let $\ps \colon A_0 \to C^* (\Z, A, \af)$
be as in the conclusion of Lemma~\ref{EmbedA}.
Suppose that $p, \, q \in \ps (A_0)$ are \pj s
such that $\ta (p) < \ta (q)$
for all tracial states $\ta$ on $C^* (\Z, A, \af).$
Then there exists a \pj\   $r \in \ps (A_0)$ such that $r \leq q$
and $r$ is \mvnt\  to $p$ in $C^* (\Z, A, \af).$
\end{lem}

\begin{proof}
If the \pj\  $f$ as in Lemma~\ref{EmbedA} is zero,
then $A_0 = A$ and $\ps = \io.$
So the statement follows from Proposition~\ref{OrdPjCrPrd}.

Otherwise, following the proof of Lemma~\ref{EmbedA},
let $e_{j, k},$ for $0 \leq j, \, k \leq n,$
be the matrix units in $M_{n + 1}.$
Also let $\io \colon A \to C^* (\Z, A, \af)$ be the inclusion,
and let $D = \io (A)$ and $D_0 = \ps (A_0).$
Since $a \in f A f$ implies $\io (a) = \ps (e_{1, 1} \otimes a),$
the algebra $E = \io (f A f)$ is a \hsa\  of both $D$ and $D_0.$

Now let $p, \, q \in D_0$ be \pj s as in the hypotheses.
Since $D_0 = \ps (A_0)$ is simple, there is $m$ such that
\[
(1, 0, \ldots, 0) \precsim
 (\io (f), \, \io (f), \, \ldots, \, \io (f))
\]
in $M_m (D_0).$
We identify $D$ and $D_0$ with corners in $M_m (D)$ and $M_m (D_0)$
in the usual way.
Then, in particular, there exist \pj s
\[
p_0, \, q_0 \leq (\io (f), \, \io (f), \, \ldots, \, \io (f))
\]
in $M_m (D_0)$ such that $p \sim p_0$ and $q \sim q_0$ in $M_m (D_0).$
Clearly $p_0, \, q_0 \in M_m (E) \subset M_m (D),$
and satisfy $\ta (p_0) < \ta (q_0)$
for $\ta \in T (C^* (\Z, A, \af)).$
Because $D = \io (A),$
Proposition~\ref{OrdPjCrPrd} provides $r_0 \in M_m (D)$
such that $p_0 \sim r_0$ in $M_m (C^* (\Z, A, \af))$ and $r_0 \leq q_0.$
Then $r_0 \in M_m (E) \subset M_m (D_0).$
Choose $s \in M_m (D_0)$ such that $s^* s = q_0$ and $s s^* = q.$
Then $r = s r_0 s^* \in M_m (D_0)$ satisfies $p \sim r$ in
$M_m (C^* (\Z, A, \af))$ and $r \leq q.$
Since $p, \, q \in C^* (\Z, A, \af),$ we in fact get
$p \sim r$ in $C^* (\Z, A, \af).$
\end{proof}

\begin{lem}\label{TComm}
Let $A$ be a \ca, let $p, \, q \in A$ be \pj s,
let $\ta$ be a tracial state on $A,$
and let $g \colon [0, 1] \to \R$ be a \cfn.
Then $\ta (g (p q p)) = \ta (g (q p q)).$
\end{lem}

\begin{proof}
The conclusion is true when $g (t) = t^n,$ since
\[
\ta ((p q p)^n) = \ta ((p q p)^{n - 1} (p q)(q p))
 = \ta ((q p) (p q p)^{n - 1} (p q)) = \ta ((q p q)^n).
\]
So it also holds for any polynomial and
therefore, by approximation, for any \cfn\  $g.$
\end{proof}

\begin{lem}\label{LBTrace}
Let $g \colon [0, 1] \to [0, 1]$ be a \cfn\  such that $g (1) = 1.$
Then for every $\ep > 0$ there exists $\dt > 0$ such that
whenever $A$ is a unital \ca, $\ta$ is a tracial state on $A,$
and $p, \, q \in A$ are \pj s such that $\ta (p) > 1 - \dt,$
then $\ta (g (q p q)) > \ta (q) - \ep$
and $\ta (g (p q p)) > \ta (q) - \ep.$
\end{lem}

\begin{proof}
We prove the result for the inequality
$\ta (g (q p q)) > \ta (q) - \ep.$
Choose $\dt_0 \in (0, 1)$ such that
$g (t) > 1 - {\textstyle{\frac{1}{2}}} \ep$
for all $t \in [1 - \dt_0, \, 1].$
Then set $\dt = {\textstyle{\frac{1}{2}}} \ep \dt_0.$
Let $A, \, \ta, \, p, \, q$ be as in the hypotheses.

We first estimate $\ta (q p q),$ as follows.
We have $\ta (q p q) + \ta (q (1 - p) q)  = \ta (q)$ and
\[
\ta (q (1 - p) q) = \ta ( (1 - p) q (1 - p) ) \leq \ta (1 - p) < \dt,
\]
so that $\ta (q p q) > \ta (q) - \dt.$

Now let $\mu$ be the measure on $X = \spec (q p q)$
corresponding to the functional on $C (X)$ defined by
$h \mapsto \ta (h (q p q)),$
with the functional calculus evaluated in $q A q.$
This measure has total mass $\ta (q).$
With $E = [1 - \dt_0, \, 1],$ we have
\begin{align*}
\ta (q) - \dt
  & < \ta (q p q)
    = \int_0^1 t \, d \mu (t)
    \leq (1 - \dt_0)
     \mu  ([0, 1] \setminus E) + \mu ( E)  \\
  & = (1 - \dt_0) [ \ta (q) - \mu ( E) ] + \mu ( E)
    = (1 - \dt_0) \ta (q) + \dt_0 \mu ( E).
\end{align*}
Rearranging this gives
\[
\mu ( E)
 > \ta (q) - \frac{\dt}{\dt_0}
 = \ta (q) - {\textstyle{\frac{1}{2}}} \ep.
\]

Since $g (t) > 1 - {\textstyle{\frac{1}{2}}} \ep$ for $t \in E,$
we now get
\begin{align*}
\ta (g (q p q) )
  & = \int_0^1 g (t) \, d \mu (t)
    \geq \left( 1 - {\textstyle{\frac{1}{2}}} \ep \right)
           \mu ( E)  \\
  & > \left( 1 - {\textstyle{\frac{1}{2}}} \ep \right)
      \left( \ta (q) - {\textstyle{\frac{1}{2}}} \ep \right)
    > \ta (q) - {\textstyle{\frac{1}{2}}} \ep
               - {\textstyle{\frac{1}{2}}} \ep \ta (q).
\end{align*}
Since $\ta (q) \leq 1,$ this gives
$\ta (g (q p q) ) > \ta (q) - \ep,$ as desired.

The result with the inequality $\ta (g (p q p)) > \ta (q) - \ep$
now follows from Lemma~\ref{TComm}.
\end{proof} % of Lemma~\ref{LBTrace}

\begin{lem}\label{AI}
Let $\dt > 0.$
Then there exists a \cfn\  $g \colon [0, 1] \to [0, 1]$
such that $g (0) = 0,$ $g (1) = 1,$
and whenever $A$ is a \ca\  with real rank zero and
$a \in A$ is a positive element with $\| a \| \leq 1,$
then there is a \pj\  $e \in {\overline{a A a}}$
such that $g (a) e  = e$ and $\| e a - a \| < \dt.$
\end{lem}

\begin{proof}
Choose $t_0$ such that
$0 < t_0 < \frac{1}{3} \dt.$
Let $g_0 \colon [0, 1] \to [0, 1]$ be a \cfn\  which vanishes on
$[0, t_0]$ and satisfies $| g_0 (t) - t | < \frac{1}{3} \dt$
for all $t \in [0, 1].$
Let $g \colon [0, 1] \to [0, 1]$ be any \cfn\  %
such that $g (0) = 0,$ $g (1) = 1,$ and $g g_0 = g_0.$

Let $A$ be a \ca\  with real rank zero and let
$a \in A$ be a positive element with $\| a \| \leq 1.$
Since $A$ has real rank zero,
there is a \pj\  $e \in {\overline{g_0 (a) A g_0 (a)}}$
such that $\| e g_0 (a) - g_0 (a) \| < \frac{1}{3} \dt.$
Since $\| a - g_0 (a) \| < \frac{1}{3} \dt,$
we get $\| e a - a \| < \dt.$
{}From $g g_0 = g_0$ we get $g (a) g_0 (a) = g_0 (a),$
whence $g (a) e = e.$
\end{proof} % of Lemma~\ref{AI}

The proof of the following theorem is adapted from the proofs
of Theorem~3.5 and Lemma~3.3 of~\cite{Ph10},
which in turn are based on Section~3 of~\cite{Pt5}.
However, the construction of the \pj\  $q_0$ in the proof is new.
It enables us to prove directly that
the order on projections over $C^* (\Z, A, \af)$ is
determined by traces,
rather than merely that the order on $K_0 (C^* (\Z, A, \af))$
is determined by traces.

\begin{thm}\label{K0OrdDetByTraces}
Let $A$ be a simple unital C*-algebra with real rank zero, and
suppose that the order on projections over $A$ is
determined by traces.
Let $\af \in \Aut (A)$ have the \tRp.
Then the order on projections over $C^* (\Z, A, \af)$ is
determined by traces.
\end{thm}

\begin{proof}
We claim that it suffices to show that if $q, \, r \in C^* (\Z, A, \af)$
are \pj s such that $\ta (q) < \ta (r)$
for all tracial states $\ta$ on $C^* (\Z, A, \af),$
then $q \precsim r.$
Indeed, it is easy to check that the action
$\id_{M_n} \otimes \af$ on $M_n \otimes A$ again has the \tRp,
so the result applies to \pj s in $M_n \otimes C^* (\Z, A, \af)$
as well,
and this version implies the statement of the theorem.

Accordingly, let $q, \, r \in C^* (\Z, A, \af)$
be \pj s such that $\ta (q) < \ta (r)$
for all tracial states $\ta$ on $C^* (\Z, A, \af).$
Since the tracial state space is weak* compact,
there is $\ep > 0$ such that
$\ta (r) - \ta (q) > \ep$ for all tracial states $\ta.$

Choose $\et > 0$ so small that
whenever $B$ is a \ca\  and $e, \, f \in B$ are \pj s
such that $\| e f - f \| < \et,$ then $f \precsim e.$
Choose \cfn s $g_1, \, g_2 \colon [0, 1] \to [0, 1]$
such that
\[
g_1 (0) = g_2 (0) = 0, \,\,\,\,\,\,
g_1 (1) = g_2 (1) = 1, \,\,\,\,\,\, g_1 g_2 = g_2,
\]
and $| g_1 (t) - t | < {\textstyle{\frac{1}{4}}} \et$
for all $t \in [0, 1].$
Choose a \cfn\  $g \colon [0, 1] \to [0, 1]$ as in Lemma~\ref{AI}
with $\frac{1}{8} \et^2$ in place of $\dt.$

Choose $\dt > 0$ so small that whenever
$B$ is a \ca\  and $a, \, b \in B$ are positive
elements with
\[
\| a \|, \, \| b \| \leq 1 \andeqn \| a - b \| < \dt,
\]
then
\[
\| g_1 (a) - g_1 (b) \| < {\textstyle{\frac{1}{4}}} \et,
\,\,\,\,\,\,
\| g_2 (a) - g_2 (b) \| < {\textstyle{\frac{1}{21}}} \ep,
\andeqn
\| g (a) - g (b) \| < {\textstyle{\frac{1}{6}}} \ep.
\]
We also require $\dt < {\textstyle{\frac{1}{2}}} \et.$

Apply Lemma~\ref{LBTrace} with $g_1$ in place of $g$
and with ${\textstyle{\frac{1}{21}}} \ep$ in place of $\ep,$
obtaining a number $\dt_0 > 0.$
Choose an integer $N \geq \max (\dt_0^{-1}, 6 \ep^{-1}).$

Apply Lemma~\ref{Basic1}
with $\{ q, \, r \}$ in place of $F,$
with ${\textstyle{\frac{1}{2}}} \dt$ in place of $\ep,$
with $N$ as given,
and with $1$ in place of $z.$
We obtain a \pj\  $e \in A \subset C^* (\Z, A, \af),$
a unital subalgebra $D \subset e C^* (\Z, A, \af) e,$
a \pj\  $p \in D \cap A,$
a \pj\  $f \in A,$ and an isomorphism
$\ph \colon M_n \otimes f A f \to D,$ satisfying the
conditions~(1) through~(6) there.

In the next several paragraphs, we construct a \pj\  $r_0 \in D$
such that $r_0 \precsim r$ and
$\ta (r_0) > \ta (r) - {\textstyle{\frac{1}{3}}} \ep$
for every tracial state $\ta$ on $C^* (\Z, A, \af).$

By the choice using Lemma~\ref{Basic1},
there exists $x \in D$ such that
$\| r p - x \| < {\textstyle{\frac{1}{2}}} \dt$ and $\| x \| \leq 1,$
so that $\| r p r - x x^* \| < \dt.$
Since $x \in D \cong M_n \otimes f A f,$ which has real rank zero,
Lemma~3.2 of~\cite{Ph10} provides a \pj\  $r_0 \in D$ such that
\[
g_1 (x x^*) r_0 = r_0
\andeqn
\| r_0 g_2 (x x^*) - g_2 (x x^*) \| < {\textstyle{\frac{1}{21}}} \ep.
\]
We claim that $r_0 \precsim r,$ and we prove this by
showing that $\| r r_0 - r_0 \| < \et.$
The choice of $\dt$ and the estimate $\| r p r - x x^* \| < \dt$
imply that
$\| g_1 ( r p r) - g_1 ( x x^*) \| < {\textstyle{\frac{1}{4}}} \et.$
So $g_1 (x x^*) r_0 = r_0$ implies
$\| g_1 ( r p r) r_0 - r_0 \| < {\textstyle{\frac{1}{4}}} \et,$
and from $| g_1 (t) - t | < {\textstyle{\frac{1}{4}}} \et$
we then get $\| r p r r_0 - r_0 \| < {\textstyle{\frac{1}{2}}} \et.$
Now
\[
\| r r_0 - r_0 \|
  \leq \| r \| \cdot \| r_0 - r p r r_0 \| + \| r^2 p r r_0 - r_0 \|
  < \et,
\]
as desired.
This proves the claim.

Now let $\ta$ be a tracial state on $C^* (\Z, A, \af).$
We obtain a lower bound on $\ta (r_0).$
The choice of $\dt$ and the estimate $\| r p r - x x^* \| < \dt$
imply that
$\| g_2 ( r p r) - g_2 ( x x^*) \| < {\textstyle{\frac{1}{21}}} \ep.$
So
$\| r_0 g_2 (x x^*) - g_2 (x x^*) \| < {\textstyle{\frac{1}{21}}} \ep$
implies
$\| r_0 g_2 ( r p r) - g_2 ( r p r) \|
               < {\textstyle{\frac{3}{21}}} \ep,$
whence
$\| r_0 g_2 ( r p r) r_0 - g_2 ( r p r) \|
                          < {\textstyle{\frac{6}{21}}} \ep.$
Therefore
\[
\ta (r_0) \geq \ta (r_0 g_2 ( r p r) r_0)
 > \ta ( g_2 ( r p r) ) - {\textstyle{\frac{6}{21}}} \ep.
\]
Now the choice using Lemma~\ref{Basic1}
implies $\ta (1 - p) \leq N^{-1} \ta (p) \leq N^{-1} < \dt_0,$
so $\ta (p) > 1 - \dt_0,$
and the choice using Lemma~\ref{LBTrace} gives
$\ta ( g_2 ( r p r) ) > \ta (r) - {\textstyle{\frac{1}{21}}} \ep.$
Thus
$\ta (r_0) > \ta (r) - {\textstyle{\frac{7}{21}}} \ep
        = \ta (r) - {\textstyle{\frac{1}{3}}} \ep.$
We have proved that $r_0$ is the required \pj.

We now construct a \pj\  $q_0 \in (1 - p) + p D p$
such that $q \precsim q_0$ and
$\ta (q_0) < \ta (q) + {\textstyle{\frac{1}{3}}} \ep$
for every tracial state $\ta$ on $C^* (\Z, A, \af).$
The method is similar to the construction of $r_0$
but is a bit more complicated.

By the choice using Lemma~\ref{Basic1},
there exists $x \in D$ such that
$\| p q - x \| < {\textstyle{\frac{1}{2}}} \dt$ and $\| x \| \leq 1.$
Replacing $x$ by $p x,$ we may assume in addition that $p x = x.$
Then $x x^* \in p D p$ and $\| p q p - x x^* \| < \dt.$
Since $D \cong M_n \otimes f A f,$ which has real rank zero,
we may apply the choice of $g$ to find a \pj\  %
$q_1 \in {\overline{x x^* D x x^*}} \subset p D p$ such that
\[
g (x x^*) q_1 = q_1
\andeqn
\| q_1 x x^* - x x^* \| < \ts{ \frac{1}{8} } \et^2.
\]
Now set $q_0 = 1 - p + q_1 \in (1 - p) + p D p.$

We estimate $\| q_0 q - q \|.$
First,
\[
\| q_1 x - x \|^2
   \leq  \| q_1 x x^* - x x^* \| \cdot\| q_1^* \|
                 + \| q_1 x x^* - x x^* \|
   \leq 2 \| q_1 x x^* - x x^* \|
   < \ts{ \frac{1}{4} } \et^2.
\]
So $\| q_1 x - x \| < \ts{ \frac{1}{2} } \et.$
Now, using $q_1 p = q_1$ at the second step,
\begin{align*}
\| q_0 q - q \|
 & = \| (1 - p) q + q_1 q - q \|
   = \| q_1 p q - p q \|  \\
 & \leq 2 \| p q - x \| + \| q_1 x - x \|
   < \dt + \ts{ \frac{1}{2} } \et
   \leq \ts{ \frac{1}{2} } \et + \ts{ \frac{1}{2} } \et
   = \et.
\end{align*}
It follows from the choice of $\et$ that $q \precsim q_0.$

Now we estimate the values of the tracial states on $q_0.$
Let $\ta$ be any tracial state on $C^* (\Z, A, \af).$
The estimate $\| p q p - x x^* \| < \dt$ and the choice of $\dt$
imply $\| g (p q p) - g (x x^*) \| < {\textstyle{\frac{1}{6}}} \ep.$
Using inequality in the \ca\  at the third step,
this estimate at the fourth step,
Lemma~\ref{TComm} at the fifth step,
and $g (q p q) \leq q$ at the sixth step, we get
\begin{align*}
\ta (q_1)
 & = \ta ( q_1 g (x x^*) q_1)
   = \ta \ts{ \left( g (x x^*)^{1/2} q_1 g (x x^*)^{1/2} \right) }
   \leq \ta (g (x x^*))  \\
 & < \ta (g (p q p)) + {\textstyle{\frac{1}{6}}} \ep
   = \ta (g (q p q)) + {\textstyle{\frac{1}{6}}} \ep
   \leq \ta (q) + {\textstyle{\frac{1}{6}}} \ep.
\end{align*}
For the same reason that we had $\ta (1 - p) < \dt_0$ in the
proof of the estimate for $\ta (r_0),$
we also have $\ta (1 - p) < {\textstyle{\frac{1}{6}}} \ep.$
Therefore
\[
\ta (q_0) = \ta (1 - p) + \ta (q_1)
 < \ta (q) + {\textstyle{\frac{1}{3}}} \ep.
\]
We have proved that $q_0$ is the required \pj.

Apply Lemma~\ref{EmbedA} with $\ph \colon M_n \otimes f A f \to D$
and the \pj\  $e$ as given.
We obtain $A_0$
and a unital \hm\  $\ps \colon A_0 \to C^* (\Z, A, \af).$
We note that $\ps (A_0)$ contains $D,$ and hence $r_0$;
also $1, \, p \in \ps (A_0)$ so
$q_0 \in (1 - p) + p D p \subset \ps (A_0).$
For every tracial state $\ta$ on $C^* (\Z, A, \af)$ we have
\[
\ta (r_0) - \ta (q_0)
 > \left( \ta (r) - {\textstyle{\frac{1}{3}}} \ep \right)
     - \left( \ta (q) + {\textstyle{\frac{1}{3}}} \ep \right)
 > {\textstyle{\frac{1}{3}}} \ep.
\]
Therefore Lemma~\ref{BasicOrder} applies, and shows that
$q_0 \precsim r_0$ in $C^* (\Z, A, \af).$
Thus,
\[
q \precsim q_0 \precsim r_0 \precsim r,
\]
as was to be proved.
\end{proof} % of Theorem~\ref{K0OrdDetByTraces}

\section{Real rank of crossed products}\label{Sec:RR}

\indent
In this section, we prove that if
$A$ is a simple unital C*-algebra with real rank zero such that
the order on projections over $A$ is
determined by traces, and if $\af \in \Aut (A)$ has the \tRp,
then $C^* (\Z, A, \af)$ has real rank zero,
and every tracial state on $C^* (\Z, A, \af)$ is induced
from an $\af$-invariant tracial state on $A.$
The methods are adapted from Section~4 of~\cite{Ph10}.
The lemma used to find a suitable \pj\  which approximately
commutes with a selfadjoint element
is considerably harder in our context.
To prove it, we start with the following lemma.
It will follow from arguments in~\cite{Cu},
where it is proved that if $c \in A$ then
${\overline{c^* A c}} \cong {\overline{c A c^*}}.$

\begin{lem}\label{chp3:lem2}
Let $A$ be a \ca, let $p \in A$ be a \pj, and let $c \in A.$
Suppose that ${\overline{c^* A c}} \subset p A p.$
Then for any projection $r \in {\overline{c A c^*}},$
we have $r \precsim p.$
\end{lem}

\begin{proof}
For each $\ep > 0$
define a \cfn\  $f_{\ep} \colon [0, \infty) \to [0, 1]$ by
\[
f_{\ep} (t) = \left\{ \begin{array}{ll}
                   0  & \hspace{3em} t \leq \frac{\ep}{2}  \\
 \frac{2}{\ep} \left( t - \frac{\ep}{2} \right)
                      & \hspace{3em} \frac{\ep}{2} \leq t \leq \ep \\
                   1  & \hspace{3em} \ep \leq t.
                   \end{array}
                   \right.
\]
Recall that $| c | = (c^* c)^{1/2},$ so that $| c^* | = (c c^*)^{1/2},$
and (analogously to~1.2 of~\cite{Cu})
\[
{\overline{c A c^*}}
 = {\overline{\sucup{\ep > 0} f_{\ep} (| c^* |) A f_{\ep} (| c^* |)}}.
\]

Now let $r \in {\overline{c A c^*}}$ be a \pj; we prove $r \precsim p.$
\Wolog\  there is $\ep > 0$ such that
\[
r \in {\overline{ f_{\ep} (| c^* |) A f_{\ep} (| c^* |) }}.
\]
Then $f_{\ep/2} (|c^* |) r = r,$ whence
\[
r \in f_{\ep/2} (|c^* |) A f_{\ep/2} (|c^* |).
\]
By~1.4 of~\cite{Cu} there is $z \in A$ such that
the map $a \mapsto z a z^*$ is an isomorphism
\[
f_{\ep/2} (| c^* |) A f_{\ep/2} (| c^* |) \to
f_{\ep/2} (|c|) A f_{\ep/2} (|c|).
\]
In particular, $z r z^*$ is a \pj\  and
$r z^*$ is a partial isometry from $r$ to $z r z^*.$
Since
\[
f_{\ep/2} (|c|) A f_{\ep/2} (|c|) \subset
{\overline{c^* A c}} \subset p A p,
\]
we have $r \sim z r z^* \leq p.$
\end{proof}

The following lemma is surely true for actions of
(not necessarily abelian) compact groups,
and may even be known.

\begin{lem}\label{InvAppId}
Let $A$ be a unital \ca, let $G$ be a finite abelian group,
let $\af \colon G \to \Aut (A)$ be an action of $G$ on $A,$
and let $A^{\af}$ be the fixed point algebra.
Then every approximate identity for $A^{\af}$
is also an approximate identity for $A.$
\end{lem}

\begin{proof}
Let $( e_{\ld} )_{\ld \in \Ld}$
be an approximate identity for $A^{\af}.$
For $\ta \in {\widehat{G}},$ let
\[
A_{\ta} = \{ a \in A \colon
  {\mbox{$\af_g (a) = \ta (g) a$ for all $g \in G$}} \}.
\]
Since $A$ is the direct sum of the subspaces $A_{\ta},$
it suffices to prove that $\lim_{\ld} a e_{\ld} = a$
for all $a \in A_{\ta}.$

So let $a \in A_{\ta}.$
Then
\[
\| ( a - a e_{\ld})^* ( a - a e_{\ld}) \|
  = \| ( a^* a - a^* a e_{\ld}) - e_{\ld} ( a^* a - a^* a e_{\ld}) \|
  \leq 2 \| a^* a - a^* a e_{\ld} \|.
\]
Since $a^* a \in A^{\af},$
we have $\lim_{\ld} \| a^* a - a^* a e_{\ld} \| = 0.$
\end{proof}

\begin{lem}\label{InnerOnHer}
Let $A$ be a \ca, let $G$ be a topological group,
and let $\af \colon G \to \Aut (A)$ be an action of $G$ on $A$
which is inner in the sense that there is
a strictly \ct\  group homomorphism $g \mapsto u (g) \in U (M (A)),$
the unitary group of the multiplier algebra of $A,$
such that $\af_g (a) = u (g) a u (g)^*$ for all $g \in G.$
Then the restriction of the action to any invariant \hsa\  $B \subset A$
is also inner in the same sense.
\end{lem}

\begin{proof}
It suffices to show that if $g \in G$ then the pair $(L, R)$ of linear
maps on $B,$ given by $L (x) = u (g) x$ and $R (x) = x u (g),$
is a multiplier of $B.$
The only nontrivial part is that
$L (B) \subset B$ and $R (B) \subset B.$
For the first, $x \in B$ implies
\[
L (x) L (x)^* = u (g) x x^* u (g)^* \in \af_g (B) = B
\]
and
\[
L (x)^* L (x) = x^* u (g)^* u (g) x = x^* x \in B,
\]
whence $L (x) \in B.$
The proof of the second is similar.
\end{proof}

\begin{lem}\label{chp3:lem3}
Let $A$ be a unital \ca\  with real rank zero.
Let $m, n, N \in \N$ satisfy $(2 n + 1) m \leq N.$
Let $a \in M_N (A)$ be a selfadjoint
element with $\| a \| \leq 1.$
Let $( e_{i, j} )_{i, j = 1}^N$ be the standard
system of matrix units of $M_N ( \C ).$
Then there is a projection $q \in M_N (A)$ such that
\[
\sum_{k = 1}^{m} e_{k, k} \otimes 1 \leq q,  \,\,\,\,\,\,
q \precsim \sum_{k = 1}^{(2 n + 1) m} e_{k, k} \otimes 1, \andeqn
\| q a - a q \| < \frac{1}{n}.
\]
\end{lem}

\begin{proof}
Set
\[
p = \sum_{k = 1}^{m} e_{k, k}
\andeqn q_0 = \sum_{k = 1}^{(2 n + 1) m} e_{k, k} \otimes 1.
\]
Then we want $q$ with
\[
p \leq q \precsim q_0 \andeqn \| q a - a q \| < \ts{ \frac{1}{n} }.
\]

Since $A$ has real rank zero, so does $M_N (A).$
Set
\[
\ep = \frac{1}{2} \left( \frac{1}{2 n} - \frac{1}{2 n + 1} \right).
\]
Because
\[
\| a \| \leq 1 \andeqn \frac{1}{2 n + 1} + \ep > \frac{1}{2 n + 1},
\]
we can approximate $a$ by a selfadjoint element in $M_N (A)$
with finite spectrum, and then perturb its spectrum, to get
a selfadjoint element $b \in M_N (A)$ such that
\[
\| a - b \| < \frac{1}{2 n + 1} + \ep \andeqn
\spec (b) \subset
\left\{ - \frac{2 n}{2 n + 1}, \, - \frac{2 n - 2}{2 n + 1},
  \, \ldots, \, \frac{2 n - 2}{2 n + 1}, \,
\frac{2 n}{2 n + 1} \right\}.
\]
Define $u = \exp (\pi i b).$
Then $u$ is a unitary in $M_N (A)$ with $u^{2 n + 1} = 1.$
Letting $\log$ be the standard branch defined on the complement
of the negative real axis,
we furthermore have $b = - \frac{i}{\pi} \log (u).$

Let $v \in M_N (A)$ be the following unitary block matrix,
whose entries are the identity matrices of the appropriate sizes:
\[
v = \left(\begin{array}{cccc}
   0 &  &  & 1_m \\
   1_m & \ddots &  & \\
   & \ddots & \ddots & \\
   &  & 1_m & 0   \end{array} \right)
 \oplus 1_{N - (2 n + 1) m}.
\]
Define
\[
c = p + u p v^* + \cdots + u^{2 n} p (v^*)^{2 n},
\]
and let $C = {\overline{c A c^*}},$ which is a \hsa\  in $M_N (A).$

For $0 \leq k, l \leq 2 n$ we have
\[
v^k p (v^*)^l = \sum_{j = 1}^m  e_{k m + j, \, l m + j}.
\]
Since $p (v^*)^l q_0 = p (v^*)^l$ for $0 \leq l \leq 2 n,$
we get $c q_0 = c,$
so that $c^* c \in q_0 M_N (A) q_0.$
It also follows that
in the computation of $c c^*$ most of the terms cancel, and one gets
\[
c c^* = p + u p u^* + \cdots + u^{2 n} p (u^*)^{2 n}.
\]
Since $u^{2 n + 1} = 1,$ it follows that $u c c^* u^* = c c^*.$

We claim that $u C u^* = C.$
Indeed, if $x \in A$ then
\[
u c x c^* u^* = (u c) x (u c)^* \leq \| x \| (u c) (u c)^*
  = \| x \| \cdot u c c^* u^* = \| x \| \cdot c c^* \in C.
\]
This shows that $u C u^* \subset C,$
and the reverse inclusion now follows from $u^{2 n + 1} = 1.$
This proves the claim.

Write $\Z_{2 n + 1} = \Z / (2 n + 1) \Z.$
The automorphism $\Ad (u)$ generates an inner action of $\Z_{2 n + 1}$
on $M_N (A),$
and by invariance and Lemma~\ref{InnerOnHer} also an inner action
$\af \colon \Z_{2 n + 1} \to \Aut (C).$
We claim that the fixed point algebra $C^{\af}$
has real rank zero.
To see this, note that the Proposition in~\cite{Rs} shows that
$C^{\af}$ is isomorphic to a \hsa\  %
of the crossed product $C^* (\Z_{2 n + 1}, \, C, \, \af).$
Since $\af$ is inner, this crossed product
is isomorphic to the direct sum of $2 n + 1$ copies of $C,$
so has real rank zero.
Therefore so does every \hsa.
This proves the claim.

Choose $\dt_1 > 0$ such that whenever $D$ is a unital \ca\  %
and an element $x \in D$ and a unitary $v \in D$ satisfy
\[
\| x \| \leq 1, \,\,\,\,\,\,
\| x v  - v x \| < 4 \dt_1, \andeqn
\spec (v)
  \subset \left\{ \exp (i \te) \colon
        - \frac{2 n}{2 n + 1}\leq \te \leq \frac{2 n}{2 n + 1} \right\},
\]
then $\| x \log (v) - \log (v) x \| < \pi \ep.$
Choose $\dt_2 > 0$ such that whenever $D$ is a unital \ca\  %
and projections $e, f \in D$ satisfy $\| e - f \| < \dt_2,$
then there exists a unitary $w \in D$
such that $\| w - 1 \| < \dt_1$ and $w e w^* = f.$
Choose $\dt_3 > 0$ such that
$\dt_3 < \frac{1}{4} \min ( \dt_2, 1 ).$

We now claim that there is a projection $e \in C$
such that $u e u^* = e$ and $\| e p - p \| < \dt_3.$
To prove this,
first observe that the formula for $c c^*$ above
shows that $p \leq c c^*,$
whence $p \in C.$
Next, use the fact that $C^{\af}$ has real rank zero
to find (Theorem~2.6 of~\cite{BP}) an approximate identity
in $C^{\af}$ consisting of \pj s.
Lemma~\ref{InvAppId} shows that any approximate identity
for $C^{\af}$ is also an approximate identity for $C.$
Since $\af (b) = u b u^*$ for $b \in C,$ the claim follows.

{}From $\| e p - p \| < \dt_3$ we get
\[
\| (e p e)^2 - e p e \|
  = \| e p e p e - e p p e \|
  \leq \| e p \| \cdot \| e p - p \| \cdot \| e \|
  < \dt_3.
\]
Since $\dt_3 < \frac{1}{4},$ a standard argument provides
a projection $f \in {\overline{(e p e) A (e p e)}}$
such that $\| e p e - f \| < 2 \dt_3.$
Then
\[
\| p - f \|
 \leq \| (p - e p)^* \| + \| p - e p \| \cdot \| e \| + \| e p e - f \|
 < 4 \dt_3 < \dt_2.
\]
Note that $e f = f,$ so $e \geq f$;
also, $e \precsim q_0$ by Lemma~\ref{chp3:lem2}.

By the choice of $\dt_2$ there is a unitary $w \in A$ such that
$\| w - 1 \| < \dt_1$ and $w f w^* = p.$
Set $q = w e w^*.$
Then $q \geq w f w^* = p.$
Since $q \sim e \precsim q_0,$ we also have $q \precsim q_0.$

It remains to show that $\| q a - a q \| < \frac{1}{n}.$
Since
\[
\| q u - u q \| = \| w e w^* u - u w e w^* \|
 \leq 4 \| w - 1 \| < 4 \dt_1,
\]
the choice of $\dt_1$ gives
\[
\| q b - b q \| = \ts{ \frac{1}{\pi} } \| q \log (u) - \log (u) q \|
  < \ep.
\]
Therefore
\[
\| q a - a q \|  \leq 2 \| a - b \| + \| q b - b q \|
 \leq 2 \left( \frac{1}{2 n + 1} + \ep \right) + \ep < \frac{1}{n}.
\]
This completes the proof.
\end{proof}

The proof of the following theorem is analogous to that of
Theorem~4.6 of~\cite{Ph10},
but is somewhat more complicated.
In~\cite{Ph10} we worked with a single ``large''
AF subalgebra; here, we only have subalgebras of real rank zero,
obtained from Lemma~\ref{Basic1},
which play the role of algebras in a direct limit decomposition
for the AF subalgebra of~\cite{Ph10}.

\begin{thm}\label{chp3:thm1}
Let $A$ be a simple unital C*-algebra with real rank zero.
Suppose that the order on projections over $A$ is
determined by traces and $\af \in \Aut (A)$ has the \tRp.
Then $C^* (\Z, A, \af)$ has real rank zero.
\end{thm}

\begin{proof}
Set $B = C^* (\Z, A, \af).$

Let $a \in B$ be selfadjoint with $\| a \| \leq 1.$
Let $\ep > 0.$
We approximate $a$ to within $\ep$
by an invertible selfadjoint element.
If $a$ is already invertible, there is nothing to prove.
Therefore we assume $0 \in \spec (a).$
Set $\ep_0 = \frac{1}{12} \ep,$
and choose a continuous function $g \colon [-1, \, 1] \to [0, 1]$
such that
\[
g (0) = 1 \andeqn
\supp (g) \subset (- \ep_0, \, \ep_0).
\]

Recalling the notation $T (B)$ from Notation~\ref{TraceNtn},
define
\[
\et = \inf_{\ta \in T (B)} \ta (g (a)).
\]
The algebra $B$ is simple by Corollary~\ref{TRPImpSimple},
so that every tracial state is faithful.
Also, $g (a)$ is a nonzero positive element,
and $T (B)$ is weak* compact.
Therefore $\et > 0.$

Use a polynomial approximation to the function $g$
to choose $\dt_{0} > 0$ such that whenever $C$ is a unital \ca\  %
and $x, y \in C_{\sa}$ satisfy $\| x \|, \| y \| \leq 2$
and $\| x - y \| < \dt_{0},$
then $\| g (x) - g (y) \| < \frac{1}{6} \et.$
Set $\dt = \min \left( \dt_0, 1, \ep_0 \right).$
Choose $n_0, N \in \N$ such that
\[
\frac{1}{N} < \frac{\et}{12} \andeqn \frac{1}{n_0} < \frac{\dt}{2}.
\]

Since $\af$ has the \tRp, we can use Lemma~\ref{Basic1}
to find a projections $e, p \in A,$
a \pj\  $f \in A,$ integers $n, \, m > 0,$
a unital C*-subalgebra $D \subset e B e,$
and an isomorphism $\ph \colon D \to M_n \otimes f A f,$
such that $p \in D,$ such that
\[
\frac{2 m}{n}
  < \min \left( \frac{1}{2 n_0 + 1}, \,
             \frac{\et}{12 (2 n_0 + 1)} \right),
\]
such that
\[
\ph (p) = \sum_{j = m + 1}^{n - m} e_{j, j} \otimes 1_{f A f}
 \in M_n \otimes f A f,
\]
such that
\[
\dist (p a, D) < \ts{ \frac{1}{2} } \dt
\andeqn \dist (a p, D) < \ts{ \frac{1}{2} } \dt,
\]
and such that there are $N$
\mops\  $f_1, f_2, \ldots, f_N \in p D p,$
each of which is \mvnt\  in $B$ to $1 - p.$

{}From the last condition, it is evident that for every $\ta \in T (B)$
we have
\[
\ta (1 - e) \leq \ta (1 - p) \leq \frac{ \ta (p) }{N} \leq \frac{1}{N}
   < \frac{\et}{12}.
\]
Moreover,
\[
\ta (e - p) = \left( \frac{2 m}{n} \right) \ta (e)
  \leq \frac{2 m}{n} < \frac{\et}{12 (2 n_0 + 1)}.
\]

Set
\[
x = a - (1 - p) a (1 - p) = p a + (1 - p) a p.
\]
Choose $x_1, x_2 \in D$ such that
\[
\| p a - x_1 \| < \ts{ \frac{1}{2} } \dt
\andeqn
\| a p - x_2 \| < \ts{ \frac{1}{2} } \dt.
\]
We arrange to replace $x_1$ and $x_2$ by a single selfadjoint element.
Since $p \in D$ and $D$ is a unital subalgebra of $e B e,$
we have
\[
(1 - p) x_2 = (1 - p) e x_2 = (e - p) x_2 \in D.
\]
So $p x_1, \, (1 - p) x_2 p \in D.$
Set $d = p x_1 + (1 - p) x_2 p \in D,$ and observe that
\[
\| d - x \|
  \leq \| p \| \cdot \| p a - x_1 \|
        + \| 1 - p \| \cdot \| a p - x_2 \| \cdot \| p \|
  < \ts{ \frac{1}{2} } \dt + \ts{ \frac{1}{2} } \dt = \dt.
\]
Then set $a_0 = a - x + \ts{ \frac{1}{2} } (d + d^*),$ which satisfies
\[
a_0^* = a_0, \,\,\,\,\,\,
a_0 - (1 - p) a_0 (1 - p) = \ts{ \frac{1}{2} } (d + d^*) \in D,
\andeqn \| a - a_0 \| < \dt.
\]

Next, we replace $p$ by a smaller \pj\  (which will be called $1 - q$)
which approximately commutes with $a_0.$
Let $z \in M_n \subset M_n \otimes f A f$
be a permutation unitary such that
\[
z \ph (e - p) z^* = \sum_{j = 1}^{2 m} e_{j, j} \otimes 1_{f A f}.
\]
Apply Lemma~\ref{chp3:lem3} with $f A f$ in place of $A,$
with $z \ph (d) z^*$ in place of $a,$
with $n$ in place of $N,$
with $2 m$ in place of $m,$
and with $n_0$ in place of $n.$
Note that
\[
(2 n_0 + 1) \cdot 2 m
  = \left( \frac{2 m}{n} \right) (2 n_0 + 1) n
  < \left( \frac{1}{2 n_0 + 1} \right) (2 n_0 + 1) n = n,
\]
as required in the hypotheses of Lemma~\ref{chp3:lem3}.
We obtain a \pj\  $q_0 \in M_n \otimes f A f$
such that $z \ph (e - p) z^* \leq q_0,$
such that $[q_0] \leq (2 n_0 + 1) [\ph (e - p)]$ in $K_0 (f A f),$
and such that
\[
\| [q_0, \, z \ph (d) z^*] \| < \ts{ \frac{1}{n_0} }
     < \ts{ \frac{1}{2} } \dt.
\]
Set $q = 1 - e + \ph^{-1} (z^* q_0 z).$
We estimate $\| [ q, a_0] \|.$
Since $z \ph (e - p) z^* \leq q,$ we get $q \geq 1 - p.$
In particular, $[q, \, (1 - p) a_0 (1 - p) ] = 0,$
whence $[q, a_0] = [q, d] = [\ph^{-1} (z^* q_0 z), \, d].$
Thus
\[
\| [q, a_0] \|
  = \| [q_0, \, z \ph (d) z^*] \|
  < \ts{ \frac{1}{2} } \dt.
\]

Let $\ta \in T (B)$; we estimate $\ta (q).$
We start with $\ta (\ph^{-1} (z^* q_0 z) ).$
We know $[\ph^{-1} (q_0)] \leq (2 n_0 + 1) [e - p]$ in $K_0 (B),$
so, using the estimate above on $\ta (e - p),$ we get
\[
\ta (\ph^{-1} (z^* q_0 z) )
  = \ta (\ph^{-1} (q_0)) \leq (2 n_0 + 1) \ta (e - p)
  < \ts{ \frac{1}{12} } \et.
\]
Using also the estimate above on $\ta (1 - e),$ we then get
\[
\ta (q) = \ta (1 - e) + \ta (\ph^{-1} (z^* q_0 z) )
 < \ts{ \frac{1}{12} } \et + \ts{ \frac{1}{12} } \et
 = \ts{ \frac{1}{6} } \et.
\]

Set $y = (1 - q) a_0 (1 - q),$
which is a selfadjoint element of $D$
with $\| y \| \leq \| a_0 \| < \| a \| + \dt \leq 2.$
Let $g (y)$ be the result of evaluating functional calculus in
$(1 - q) D (1 - q).$
Since $D$ has real rank zero,
there is a \pj\  $r \in {\overline{g (y) D g (y)}}$
such that $\| r g (y) - g (y) \| < \ts{ \frac{1}{6} } \et.$

Let $\ta \in T (B)$;
we claim that $\ta (r) > \ta (q).$
First, $\| r g (y) r - g (y) \| < \ts{ \frac{1}{3} } \et.$
Since $g \leq 1$ we get $r g (y) r \leq r,$
so that
\[
\ta (r) \geq \ta (r g (y) r)
          > \ta (g (y)) - \ts{ \frac{1}{3} } \et.
\]
Next,
\[
\| a_0 - (q a_0 q + y) \|
  \leq \| q a_0 (1 - q) \| + \| (1 - q) a_0 q \|
  \leq 2 \| [q, a_0 ] \| < 2 \cdot \ts{ \frac{1}{2} } \dt
  = \dt \leq \dt_0.
\]
Let $g (q a_0 q)$ be the result of evaluating functional calculus in
$q B q.$
Then orthogonality of $q a_0 q$ and $y,$
together with the choice of $\dt_0,$ gives
\[
\| g (a_0) - [g (q a_0 q) + g (y) ] \|
  = \| g (a_0) - g (q a_0 q + y)  \|
  < \ts{ \frac{1}{6} } \et.
\]
Since $g (q a_0 q) \leq q,$
the estimate $\ta (q) < \ts{ \frac{1}{6} } \et$ implies
\[
\ta (g (y))
 > \ta (g (a_0)) - \ta (g (q a_0 q)) - \ts{ \frac{1}{6} } \et
 \geq \ta (g (a_0)) - \ta (q) - \ts{ \frac{1}{6} } \et
 > \ta (g (a_0)) - \ts{ \frac{1}{3} } \et.
\]
Moreover, $\| a - a_0 \| < \dt \leq \dt_0$ gives
$\| g (a) - g (a_0) \| < \ts{ \frac{1}{6} } \et,$
whence $\ta (g (a_0) ) > \ta (g (a)) - \ts{ \frac{1}{6} } \et.$
By the choice of $\et$ we have $\ta (g (a)) \geq \et.$
Putting everything together, we get
\[
\ta (r) > \ta (g (y)) - \ts{ \frac{1}{3} } \et
          > \ta (g (a_0)) - \ts{ \frac{2}{3} } \et
          > \ta (g (a)) - \ts{ \frac{5}{6} } \et
          \geq \ts{ \frac{1}{6} } \et > \ta (q).
\]
This proves the claim.

Since $r \in {\overline{ g (y) B g (y) }},$
the condition on $\supp (g)$ and Lemma~4.5 of~\cite{Ph10} give
\[
\| r y - y r \| < 2 \ep_0
\andeqn \| r y r \| < \ep_0.
\]
Since $r \leq 1 - q$ and $y = (1 - q) a_0 (1 - q),$
we have $r a_0 r = r y r,$ whence $\| r a_0 r \| < \ep_0.$
Also,
\begin{align*}
\| [r, a_0] \|
& = \| r y + r (1 - q) a_0 q - [y r + q a_0 (1 - q) r] \|  \\
& \leq \| [r, y] \| + 2 \| [a_0, q ] \|
  < 2 \ep_0 + 2 \left( \ts{ \frac{1}{2} } \dt \right)
  \leq 3 \ep_0.
\end{align*}

Define
\[
a_1 = (1 - q - r) a_0 (1 - q - r) + q a_0 q.
\]
We estimate $\| a_1 - a \|.$
We have
\begin{align*}
a_0 - a_1
 & = (1 - q - r) a_0 q + q a_0 (1 - q - r)
       + (1 - q - r) a_0 r + r a_0 (1 - q - r)    \\
 & \hspace*{2em} \mbox{} + r a_0 r + r a_0 q + q a_0 r.
\end{align*}
So, using
$\| [q, a_0] \| \leq \ts{ \frac{1}{2} } \dt
                \leq \ts{ \frac{1}{2} } \ep_0,$
we get
\[
\| a_0 - a_1 \|
  \leq 4 \| [q, a_0] \| + 2 \| [r, a_0] \| + \| r a_0 r \|
  < 2 \ep_0 + 6 \ep_0 + \ep_0 = 9 \ep_0.
\]
Since $\| a_0 - a \| < \dt \leq \ep_0,$ it follows that
\[
\| a_1 - a \| < 10 \ep_0.
\]

Let $A_0$ and $\ps \colon A_0 \to C^* (\Z, A, \af)$ be as in
Lemma~\ref{EmbedA}, using $\ph^{-1}$ in place of $\ph$
and with $e$ as above.
Then $\ph^{-1} (z^* q_0 z) \in D \subset \ps (A_0)$
and $1 - e \in \ps (A_0),$
so $q \in \ps (A_0)$;
also $r \in D \subset \ps (A_0)$ by construction.
We proved above that $\ta (r) > \ta (q)$ for all $\ta \in T (B).$
So Lemma~\ref{BasicOrder} implies $q \precsim r$ in $B.$
Therefore Lemma~8 of~\cite{Gd1} provides
an invertible selfadjoint element
$b_1 \in (q + r) B (q + r)$ such that
$\| b_1 - q a_0 q \| < \ep_0.$
Also, by construction, we have
\[
1 - q, \, r, \,
y = (1 - q) a_0 (1 - q) \in D,
\]
so $(1 - q - r) a_0 (1 - q - r) \in D.$
Since $D$ has real rank zero, there is an invertible selfadjoint element
$b_2 \in D$ such that
\[
\| b_2 - (1 - q - r) a_0 (1 - q - r) \| < \ep_0.
\]
Then $b_1 + b_2$ is an invertible selfadjoint element of $B,$
and satisfies
\begin{align*}
\| (b_1 + b_2) - a \|
 & \leq \| a - a_1 \| + \| b_1 - q a_0 q \|
           + \| b_2 - (1 - q - r) a_0 (1 - q - r) \|   \\
 & < 10 \ep + \ep_0 + \ep_0 = 12 \ep_0 = \ep.
\end{align*}
This completes the proof.
\end{proof}

\begin{cor}\label{Traces}
Let $A$ be a simple unital C*-algebra with real rank zero,
and suppose that the order on projections over $A$ is
determined by traces.
Let $\af \in \Aut (A)$ have the \tRp.
Then the restriction map is a bijection from the
tracial states of $C^* (\Z, A, \af)$ to the $\af$-invariant
tracial states of $A.$
\end{cor}

\begin{proof}
Since $C^* (\Z, A, \af)$ has real rank zero by Theorem~\ref{chp3:thm1},
this follows from Proposition~2.2 of~\cite{Ks3}.
\end{proof}

\begin{cor}\label{PAlg}
Let $A$ be a simple separable nuclear unital C*-algebra
with tracial rank zero and
satisfying the Universal Coefficient Theorem.
Let $\af \in \Aut (A)$ have the \tRp.
Then $C^* (\Z, A, \af)$ satisfies
the local approximation property of Popa~\cite{Pp}
(is a Popa algebra in the sense of Definition~1.2 of~\cite{Bn}).
\end{cor}

\begin{proof}
By Corollary~5.7 and Theorem~6.8 of~\cite{LnTTR},
the order on projections over $A$ is determined by traces,
and by Theorem~3.4 of~\cite{LnTAF}, the algebra $A$ has real rank zero.
So $C^* (\Z, A, \af)$ has real rank zero by
Theorem~\ref{chp3:thm1}.
It embeds in an AF~algebra by Corollary~1 at the end of Section~3
of~\cite{Mt}, and is hence quasidiagonal.
That $C^* (\Z, A, \af)$ satisfies the local approximation property of Popa
now follows from Theorem~1.2 of~\cite{Pp}.
\end{proof}

\section{Stable rank of crossed products}\label{Sec:tsr}

\indent
In this section, we prove that if
$A$ is a simple unital C*-algebra with real rank zero and
stable rank one, such that
the order on projections over $A$ is
determined by traces, and if $\af \in \Aut (A)$ has the \tRp,
then $C^* (\Z, A, \af)$ has stable rank one.
The methods are adapted from Section~5 of~\cite{Ph10}.

\begin{lem}\label{chp4:lem1}
Let $\dt > 0.$
Then there exists a \cfn\  $g \colon [0, 1] \to [0, 1]$
such that $g (0) = 0,$ $g (1) = 1,$
and whenever $A$ is a \ca\  with real rank zero and
$a \in A$ is a positive element with $\| a \| \leq 1,$
then there is a \pj\  $e \in {\overline{a A a}}$
such that $\| e g (a) - g (a) \| < \dt$ and $\| a e - e \| < \dt.$
\end{lem}

\begin{proof}
Choose $t_0,$ $t_1$ such that
$1 - \dt < t_0 < t_1 < 1.$
Let $g \colon [0, 1] \to [0, 1]$ be any \cfn\  which vanishes on
$[0, t_1]$ and satisfies $g (1) = 1.$

Let $A$ be a \ca\  with real rank zero and let
$a \in A$ be a positive element with $\| a \| \leq 1.$
Choose a \cfn\  $h \colon [0, 1] \to [0, 1]$
which vanishes on $[0, t_0]$ and satisfies $h (t) = 1$
for $t \in [t_1, 1].$
Since $A$ has real rank zero,
there is a \pj\  $e \in {\overline{g (a) A g (a)}}$
such that $\| e g (a) - g (a) \| < \dt.$
Moreover, from $h g = g$ we get $h (a) g (a) = g (a),$
whence $h (a) e = e.$
We also have $\| a h (a) - h (a) \| < \dt$
because $| t - 1 | \leq 1 - t_0 < \dt$ whenever $h (t) \neq 0.$
Accordingly,
\[
\| a e - e \| = \| a h (a) e - h (a) e \|
  \leq \| a h (a) - h (a) \| \cdot \| e \| < \dt,
\]
as was to be proved.
\end{proof}

\begin{lem}\label{LemmaForTSR}
Let $A$ be a simple C*-algebra with real rank zero and such that
the order on projections over $A$ is determined by traces.
Let $\af \in \Aut (A)$ have the \tRp.
Let $q_1, \ldots, q_n \in C^* (\Z, A, \af)$ be nonzero \pj s,
let $a_1, \ldots, a_m \in C^* (\Z, A, \af)$ be arbitrary,
and let $\ep > 0.$
Then there exists a unital subalgebra $A_0 \subset C^* (\Z, A, \af)$
which is stably isomorphic to $A,$
a \pj\  $p \in A_0,$ nonzero \pj s $r_1, \ldots, r_n \in p A_0 p,$
and elements $b_1, \ldots, b_m \in C^* (\Z, A, \af),$
such that:
\begin{itemize}
\item[(1)]
$\| q_k r_k - r_k \| < \ep$ for $1 \leq k \leq n.$
\item[(2)]
For $1 \leq k \leq n$ there is a \pj\  $g_k \in r_k A_0 r_k$
such that $1 - p \sim g_k$ in $C^* (\Z, A, \af).$
\item[(3)]
$\| a_j - b_j \| < \ep$ for $1 \leq j \leq m.$
\item[(4)]
$p b_j p \in p A_0 p$ for $1 \leq j \leq m.$
\end{itemize}
\end{lem}

\begin{proof}
Set $B = C^* (\Z, A, \af).$

Let
\[
\et = \min_{1 \leq k \leq n}
   \left( \inf_{\ta \in T (B)} \ta (q_k) \right)
    > 0
\andeqn
\ep_0 = \min \left( \frac{\et}{5}, \, \frac{\ep^2}{2} \right).
\]
Apply Lemma~\ref{chp4:lem1} with $\ep_0$ in place of $\dt,$
obtaining a \cfn\  $g \colon [0, 1] \to [0, 1].$
Apply Lemma~\ref{LBTrace} with this function $g$ and with
$\ep_0$ in place of $\ep,$
obtaining a number $\dt > 0$ such that whenever
$\ta$ is a tracial state on $B$
and $p, \, q \in B$ are \pj s such that $\ta (q) > 1 - \dt,$
then $\ta (g (q p q)) > \ta (p) - \ep_0.$
Further choose $\ep_1 > 0$
with $\ep_1 \leq \min \left( \ep_0, \ep \right)$
and so small that whenever
$x, y \in B$ are positive elements with $\| x \|, \, \| y \| \leq 1$
and $\| x - y \| < \ep_1,$ then
$\| g (x) - g (y) \| < \ep_0.$
Then choose $\ep_2 > 0$ with $\ep_2 \leq \ep_1$ and
so small that if $x, y \in B$ are selfadjoint elements
with $\| x \|, \, \| y \| \leq 1$
and $\| x - y \| < \ep_2,$ then the positive parts $x_+$ and $y_+$
satisfy $\| x_+ - y_+ \| < \ep_1.$
Apply Lemma~\ref{Basic1} with
$F = \{ q_1, \ldots, q_n, a_1, \ldots, a_m \},$
with $\ep_2$ in place of $\ep,$
with an integer $N$ so large that $1 / N < \min (\dt, \ep_0),$
and with $z = 1.$
We obtain \pj s $e \in A \subset B$ and $f \in A,$
a unital subalgebra $D \subset e B e,$
an isomorphism $\ph  \colon M_n \otimes f A f \to D,$
and a \pj\  $p \in D,$
such that, in particular, there exist
\[
x_1, \ldots, x_n, c_1, \ldots, c_m \in D
\]
with
$\| p a_j - c_j \| < \ep_2$ for $1 \leq j \leq m,$ and
\[
\| x_k \| \leq 1 \andeqn \| p q_k - x_k \| < \ep_2
\]
for $1 \leq k \leq n.$
Moreover, $\ta (1 - p) < \min (\dt, \ep_0)$ for every $\ta \in T (B).$

Apply Lemma~\ref{EmbedA} with $\ph \colon M_n \otimes f A f \to D$
and the \pj\  $e$ as given.
We obtain a \ca\  $A_0$ which is stably isomorphic to $A$
and a unital \hm\  $\ps \colon A_0 \to C^* (\Z, A, \af).$
The subalgebra $\ps (A_0)$ will be the algebra $A_0$ called for in
the statement of the lemma.
We note that $\ps (A_0)$ contains $D,$ and hence $p.$

For $1 \leq j \leq m,$ set $b_j = a_j + p (c_j - a_j) p,$
which satisfies
\[
\| b_j - a_j \| < \ep_2 \leq \ep_1 \leq \ep
\andeqn p b_j p = p c_j p \in D \subset \ps (A_0).
\]
These are Parts~(3) and~(4) of the conclusion.

Next, for $1 \leq k \leq n,$ observe that
$\ts{\frac{1}{2}} ( p x_k p + p x_k^* p)$
is a selfadjoint element of $p D p$ of norm at most one such that
\[
\left\| p q_k p - \ts{\frac{1}{2}} ( p x_k p + p x_k^* p) \right\|
  < \ep_2.
\]
So
\[
y_k = \ts{\frac{1}{2}} ( p x_k p + p x_k^* p)_+
\]
is a positive element of $p D p$ of norm at most one such that
$\left\| p q_k p - y_k \right\| < \ep_1.$

By the choice of $g$ using Lemma~\ref{chp4:lem1},
there exists a \pj\  $r_k \in p D p \subset \ps (A_0)$ such that
\[
\| r_k y_k - r_k \| < \ep_0 \andeqn
\| r_k g (y_k) - g (y_k) \| < \ep_0.
\]
Using $r_k \leq p$ at the second step, we now have
\[
(r_k q_k - r_k) (q_k r_k - r_k)
   = r_k - r_k q_k r_k
   = r_k - r_k p q_k p r_k.
\]
Therefore
\begin{align*}
\| r_k q_k - r_k \|^2
  & = \| r_k - r_k p q_k p r_k \|
    \leq \| r_k - r_k y_k \| \cdot \| r_k \|
           + \| y_k - p q_k p \|              \\
  & < \ep_1 + \ep_0 \leq \ep^2,
\end{align*}
so $\| r_k q_k - r_k \| < \ep,$
which is Part~(1) of the conclusion.

We now estimate the traces on $r_k.$
For every $\ta \in T (B),$
we have $\ta (r_k) \geq \ta (r_k g (y_k) r_k).$
By construction we have $\| r_k g (y_k) - g (y_k) \| < \ep_0,$
whence $\| r_k g (y_k) r_k - g (y_k) \| < 2 \ep_0.$
{}From $\| y_k - p q_k p \| < \ep_1$ and the choice of $\ep_1,$
we get $\| g (y_k) - g (p q_k p) \| < \ep_0.$
Since $\ta (p) > 1 - \dt,$ the choice of $\dt$ using
Lemma~\ref{LBTrace} implies that
$\ta (g (p q_k p)) > \ta (q_k) - \ep_0.$
Combining all these, we get $\ta (r_k) > \ta (q_k) - 4 \ep_0.$
On the other hand, $\ta (1 - p) \leq \ep_0.$
Since $\ep_0 \leq \frac{1}{5} \et \leq \frac{1}{5} \ta (q_k),$
we get $\ta (r_k) > \ta (1 - p).$
Since $\ta \in T (B)$ is arbitrary,
and since $1 - p$ and $r_k$ are in $\ps (A_0),$
Lemma~\ref{BasicOrder} gives Part~(2) of the conclusion.
\end{proof} % of Lemma~\ref{LemmaForTSR}

\begin{thm}\label{chp4:thm1}
Let $A$ be a simple C*-algebra with real rank zero
and stable rank one, and such that
the order on projections over $A$ is determined by traces.
Let $\af \in \Aut (A)$ have the \tRp.
Then $C^* (\Z, A, \af)$ has stable rank one.
\end{thm}

\begin{proof}
Let $B = C^* (\Z, A, \af).$

We are going to show that every two sided zero divisor
in $B$ is a limit of invertible elements.
That is, if $a \in B$ and there are
nonzero $x, y \in B$ such that
$x a = a y = 0,$ then we show that for every $\ep > 0$
there is an invertible element $c \in B$
such that $\| a - c \| < \ep.$
Because $B$ has a faithful tracial state,
every one sided invertible element is invertible.
Therefore Theorem 3.3(a) of~\cite{Ro} will imply that
any element is a limit of invertible elements, that is,
$B$ has stable rank one.

So let $a \in B,$
let $x, y \in B$ be nonzero elements such that $x a = a y = 0,$
and let $\ep > 0.$
\Wolog\  $\| a \| \leq \frac{1}{2}$ and $\ep \leq 1.$
Since $B$ has real rank zero by Theorem~\ref{chp3:thm1},
there are are nonzero projections
\[
e \in {\overline{x^* B x}} \andeqn f \in {\overline{y B y^*}},
\]
and we have $e a = a f = 0.$
Apply Lemma~\ref{LemmaForTSR} to the nonzero \pj s $e$ and $f$
and the element $a,$ with $\ts{ \frac{1}{13} } \ep$ in place of $\ep.$
Call the resulting subalgebra $A_0,$
the resulting \pj\  $p_0,$ the resulting nonzero \pj s $e_0$ and $f_0,$
and the resulting element $x_0.$
Thus
\[
e_0, \, f_0, \, p_0 x_0 p_0 \in p_0 A_0 p_0, \,\,\,\,\,\,
1 - p_0 \precsim e_0, f_0,
\]
and
\[
\| e e_0 - e_0 \|, \, \| f f_0 - f_0 \|, \, \| a - x_0 \|
        < \ts{ \frac{1}{13} } \ep.
\]

Define $a_0 = (1 - e_0) x_0 (1 - f_0).$
We clearly have $e_0 a_0 = a_0 f_0 = 0,$
and we claim that $\| a - a_0 \| < \ts{ \frac{5}{13} } \ep.$
First, using
\[
\| a \| \leq 1
\andeqn
\| e_0 e - e_0 \| = \| e e_0 - e_0 \|  < \ts{ \frac{1}{13} } \ep,
\]
we have
\[
\| e_0 x_0 \|
 \leq \| e_0 \| \cdot \| x_0 - a \| + \| e_0 - e_0 e \| \cdot \| a \|
              + \| e_0 e a \|
 < \ts{ \frac{1}{13} } \ep + \ts{ \frac{1}{13} } \ep + 0
 = \ts{ \frac{2}{13} } \ep.
\]
Similarly, $\| x_0 f_0 \| < \ts{ \frac{2}{13} } \ep.$
Therefore
\begin{align*}
\| a - a_0 \|
 & \leq \| a - x_0 \| + \| x_0 - (1 - e_0) x_0 (1 - f_0) \|        \\
 & \leq \| a - x_0 \| + \| e_0 x_0 \|
              + \| 1 - e_0 \| \cdot \| x_0 f_0 \|
   < \ts{ \frac{1}{13} } \ep + \ts{ \frac{2}{13} } \ep
               + \ts{ \frac{2}{13} } \ep
   = \ts{ \frac{5}{13} } \ep.
\end{align*}
This proves the claim.
{}From $\| a \| \leq \frac{1}{2}$
and $\ep \leq 1$ we now get $\| a_0 \| \leq 1.$

Since $A$ has real rank zero and $A_0$ is stably isomorphic to $A,$
the algebra $A_0$ also has real rank zero.
So Proposition~1.8 of~\cite{Cu} and Lemma~\ref{chp3:lem2}
show that there is a nonzero \pj\  $r \leq e_0$ such that
$r \precsim f_0.$
Similarly, $A_0$ has stable rank one,
so in fact there is a unitary $v \in A_0$ such that $v^* r v \leq f_0.$
Then $r (a_0 v^*) = (a_0 v^*) r = 0.$

Apply Lemma~\ref{LemmaForTSR} to the nonzero \pj\  $r$
and the element $a_0 v^*,$ with $\ts{ \frac{1}{13} } \ep$
in place of $\ep.$
Call the resulting subalgebra $A_1,$
the resulting \pj\  $p_1,$ the resulting nonzero \pj\  $e_1,$
and the resulting element $x_1.$
Thus
\[
e_1, \, p_1 x_1 p_1 \in p_1 A_1 p_1, \,\,\,\,\,\,
\| r e_1 - e_1 \|, \, \| a_0 v^* - x_1 \| < \ts{ \frac{1}{13} } \ep,
\andeqn
1 - p_1 \precsim e_1.
\]

Define $a_1 = (1 - e_1) x_1 (1 - e_1).$
We clearly have $e_1 a_1 = a_1 e_1 = 0.$
Also, $p_1 a_1 p_1 = (1 - e_1) p_1 x_1 p_1 (1 - e_1) \in p_1 A_1 p_1.$
Furthermore, since still $\| a_0 v^* \| \leq 1,$ the
argument used above to prove $\| a - a_0 \| < \ts{ \frac{5}{13} } \ep$
now shows that $\| a_0 v^* - a_1 \| < \ts{ \frac{5}{13} } \ep.$
So $\| a v^* - a_1 \| < \ts{ \frac{10}{13} } \ep.$
The conclusion of Lemma~\ref{LemmaForTSR} provides $s \in B$ such that
\[
s^* s = 1 - p_1, \,\,\,\,\,\, s s^* \leq e_1, \andeqn s s^* \in A_1.
\]
Set $e_2 = s s^*$ and $w = s + s^* + p_1 - e_2.$
Since $e_2 \leq e_1 \leq p_1,$ it follows that
$w$ is a unitary satisfying
\[
w e_2 w^* = 1 - p_1, \,\,\,\,\,\, w (1 - p_1) w^* = e_2,
\andeqn w (p_1 - e_2) = p_1 - e_2.
\]

We now have
$e_2 a_1 w = 0$ and $a_1 w (1 - p_1) = a_1 e_2 w = 0.$
Therefore,
with respect to the decomposition of the identity
\[
1 = e_2 \oplus (p_1 - e_2) \oplus (1 - p_1),
\]
and with $c = (p_1 - e_2) a_1 w (p_1 - e_2)$
and suitable $x, y, z \in B,$
the element $a_1 w$ has the block matrix form
\[
a_1 w = \left(
\begin{array}{ccc}
0 & 0   & 0 \\ x & c & 0 \\ y & z   & 0
\end{array}
\right).
\]

Now use $w (p_1 - e_2) = p_1 - e_2$ and $e_2 \leq p_1$ to rewrite
\[
c = (p_1 - e_2) p_1 a_1 p_1 (p_1 - e_2)
 \in (p_1 - e_2) A_1 (p_1 - e_2).
\]
Since $(p_1 - e_2) A_1 (p_1 - e_2)$ has stable rank one,
there exists an invertible element $d \in (p_1 - e_2) A_1 (p_1 - e_2)$
such that $\| c - d \| < {\textstyle{ \frac{1}{13} }} \ep.$
Then
\[
a_2 = \left( \begin{array}{ccc}
{\textstyle{ \frac{1}{13} }} \ep e_2  & 0 & 0    \\
x      & d & 0    \\
y      & z & {\textstyle{ \frac{1}{13} }} \ep (1 - p_1)
\end{array} \right)
\]
is an invertible element in $B,$ which satisfies
$\| a_2 - a_1 w \| < {\textstyle{ \frac{3}{13} }} \ep.$
So also $a_2 w^* v$ is an invertible element in $B,$ and satisfies
\[
\| a_2 w^* v - a \|
  \leq \| a_2 - a_1 w \| + \| a_1 - a v^* \|
  < \ts{ \frac{3}{13} } \ep + \ts{ \frac{10}{13} } \ep = \ep.
\]
This is the required approximation by an invertible element.
\end{proof} % of Theorem~\ref{chp4:thm1}.

\begin{cor}\label{Canc}
Let $A$ be a simple C*-algebra with real rank zero
and stable rank one, and such that
the order on projections over $A$ is determined by traces.
Let $\af \in \Aut (A)$ have the \tRp.
Then the \pj s in $\Mi (C^* (\Z, A, \af))$ satisfy cancellation:
if $e, f, p \in \Mi (C^* (\Z, A, \af))$ are \pj s such that
$e \oplus p$ is \mvnt\  to $f \oplus p,$ then $e$ is \mvnt\  to $f.$
\end{cor}

\begin{proof}
This follows from the fact that $C^* (\Z, A, \af)$ has
stable rank one (Theorem~\ref{chp4:thm1}),
using Proposition 6.5.1 of~\cite{Bl1}.
\end{proof}

\section{Examples}\label{Sec:NonFurstEx}

\indent
In this section we give some examples of crossed products by
automorphisms with the \tRp.
The examples we are most interested in require a longer
treatment, and will appear separately~\cite{OP2}.

We believe that if an action $\af$ of $\Z$ on a simple \ca\  $A$
has the \tRp, and if $A$ has tracial rank zero,
then $C^* (\Z, A, \af)$ should again have tracial rank zero.
This would in particular imply that the crossed products by
the Furstenberg transformations on irrational rotation algebras
that we consider in~\cite{OP2} have tracial rank zero,
and also that the crossed product in Example~\ref{RPCrPrdNotAT}
has tracial rank zero.
However, we give here some examples to which such a theorem can't apply,
because neither the original algebra nor the crossed
product has tracial rank zero.

For easy reference, we state the following two results.

\begin{prp}\label{ExistRP}
There exists an automorphism $\bt$ of the $2^{\infty}$ UHF algebra $B$
which generates an action of $\Z$ with the Rokhlin property
and which is the identity on K-theory.
\end{prp}

\begin{proof}
This is implicitly proved in Sections~4 and~5 of~\cite{BKRS},
although the Rokhlin property is not explicitly mentioned there.
(See~\cite{Ks2} for an explicit proof
for the $n^{\infty}$ UHF algebra for arbitrary $n.$
Note that every automorphism of $B$ is the identity on K-theory.)
\end{proof}

\begin{prp}\label{TensWithRP}
Let $A$ be a unital \ca, and let $\af \in \Aut (A)$ be arbitrary.
Let $B$ be a unital \ca, and let $\bt \in \Aut (B)$ generate
an action of $\Z$ with the Rokhlin property.
Then $\af \otimes \bt$ generates
an action of $\Z$ on $A \otimes_{\min} B$
with the Rokhlin property.
\end{prp}

\begin{proof}
Using density of the algebraic tensor product,
one sees that it suffices to simply tensor appropriate systems of
Rokhlin \pj s for $\bt$ with $1_A.$
\end{proof}

Of course, the proof works for any tensor product on which
$\af \otimes \bt$ extends to an automorphism,
in particular for $A \otimes_{\max} B.$
The situation for the \tRp\  is much less clear.

The following example shows that
the implication (1)~implies~(5) of Theorem~6.4 of~\cite{Ks4}
is no longer valid when the action is not approximately inner.

\begin{exa}\label{RPCrPrdNotAT}
We sketch an example of an automorphism $\af$ of a simple unital
AF~algebra $A$ which has the Rokhlin property but such that
$C^* (\Z, A, \af)$ is not an AT algebra.

Let $A$ be the simple unital AF~algebra such that
$K_0 (A) \cong
 \Z \left[ \frac{1}{2} \right] \oplus \Z \left[ \frac{1}{2} \right],$
with the strict order from the first coordinate,
and with $[1] \mapsto (1, 0).$
(One checks that this is in fact a Riesz group.
See Section~7.6 of~\cite{Bl1}.)
For any $d \in \Z,$ the matrix
\[
\left( \begin{array}{cc}
  1     &  0        \\
  d     &  1
\end{array} \right)
\]
defines an automorphism of $K_0 (A)$ as a scaled ordered group.
Let $\af_0 \in \Aut (A)$ induce this automorphism on K-theory.

This automorphism need not have the Rokhlin property.
Let $B$ be the $2^{\infty}$ UHF algebra and
let $\bt \in \Aut (B)$ be as in Proposition~\ref{ExistRP}.
K-theory computations show that $A \otimes B \cong A,$
and $\af = \af_0 \otimes \bt$ does have the Rokhlin property
(by Proposition~\ref{TensWithRP}),
and induces the same map on K-theory.
Since $K_0 (A)$ has a unique state,
$A$ has a unique tracial state,
so Proposition~\ref{SRPImpTRP} shows that $\af$ has the \tRp.

The Pimsner-Voiculescu exact sequence~\cite{PV} shows that
$K_0 (C^* (\Z, A, \af))$ is isomorphic to the cokernel of the
map on $K_0 (A)$ induced by
\[
\id - \af_* = \left( \begin{array}{cc}
  0     &  0        \\
 - d    &  0
\end{array} \right).
\]
If, say, $d = 3,$ then this cokernel has torsion.
Therefore $C^*(\Z, A, \af))$ is not an AT~algebra.

On the other hand, Theorem~\ref{RPImpTRP} implies that
$\af$ generates an action with the \tRp.
So Theorem~\ref{chp3:thm1}, Theorem~\ref{chp4:thm1},
and Theorem~\ref{K0OrdDetByTraces} show that $C^* (\Z, A, \af)$
has real rank zero and stable rank one,
and that the order on \pj s over this algebra is determined by traces.
\end{exa}

The remaining examples are all on \ca s which do not have
tracial rank zero.

\begin{exa}\label{NonTAF1}
Let $n \in \{ 2, 3, \ldots, \infty \},$
let $F_n$ be the free group on $n$ generators,
and let $\af$ be any automorphism of $C^*_{\mathrm{r}} (F_n).$
(An example which is particularly interesting in this context
is to take $n = \infty$ and to take $\af$ to be induced by
an infinite order permutation of the free generators of $F_n.$
Another possibility is to have $\af$ multiply the $k$-th generating
unitary by an irrational number $\ld_k.$)
Let $B$ be the $2^{\infty}$ UHF algebra and
let $\bt \in \Aut (B)$ be as in Proposition~\ref{ExistRP}.
Then
$\af \otimes \bt$ generates an action with the Rokhlin property
by Proposition~\ref{TensWithRP}.
Since $C^*_{\mathrm{r}} (F_n)$ has a unique tracial state,
it follows from Corollary~6.6 of~\cite{Ro} that
$C^*_{\mathrm{r}} (F_n) \otimes B$ has stable rank one.
Moreover, $C^*_{\mathrm{r}} (F_n) \otimes B$ is exact,
so every quasitrace is a trace
(\cite{Hg}), whence Theorem~7.2 of~\cite{Ro2}
implies that $C^*_{\mathrm{r}} (F_n) \otimes B$ has real rank zero
and Theorem 5.2(b) of~\cite{Ro2} implies that the order on \pj s
over $C^*_{\mathrm{r}} (F_n) \otimes B$ is determined by traces.
(In fact, $K_0 (C^*_{\mathrm{r}} (F_n) \otimes B)$
is $\Z \left[ \frac{1}{2} \right]$ with its usual order.)
We can now use Proposition~\ref{SRPImpTRP} to conclude that
$\af \otimes \bt$ generates an action with the \tRp.
On the other hand, the corollary to Theorem~A1 of~\cite{Rs2}
shows that $C^*_{\mathrm{r}} (F_n)$ is not quasidiagonal,
so $C^*_{\mathrm{r}} (F_n) \otimes B$ is not quasidiagonal either.
Theorem~3.4 of~\cite{LnTAF} therefore shows that
$C^*_{\mathrm{r}} (F_n) \otimes B$ does not have tracial rank zero.
Theorem~\ref{chp3:thm1}, Theorem~\ref{chp4:thm1},
and Theorem~\ref{K0OrdDetByTraces} show that the crossed product
$C^* (\Z, \, C^*_{\mathrm{r}} (F_n) \otimes B, \, \af \otimes \bt)$
has real rank zero and stable rank one,
and that the order on \pj s over this algebra is determined by traces.
However, it does not have tracial rank zero because it
contains the nonquasidiagonal \ca\  $C^*_{\mathrm{r}} (F_n).$
\end{exa}

\begin{exa}\label{NonTAF2}
Let $A$ be the simple separable \ca\  of Theorem~7.20 of~\cite{Bn}.
This algebra has real rank zero and stable rank one,
and the order on \pj s over $A$ is determined by traces.
It also has a number of other nice properties:
it is exact,
it satisfies the Universal Coefficient Theorem,
it is approximately divisible in the sense of~\cite{BKR},
it is a direct limit of residually finite dimensional \ca s,
and it satisfies the local approximation property of Popa
(is a Popa algebra in the sense of Definition~1.2 of~\cite{Bn})
and is hence quasidiagonal (by Theorem~1.2 of~\cite{Pp}).
According to Corollary~7.21 of~\cite{Bn},
the algebra $A$ does not have tracial rank zero.

{}From the construction in the proof of Theorem~7.20 of~\cite{Bn},
we see that $A$ can be chosen to be a tensor product of some other
\ca\  with an arbitrary UHF algebra.
In particular, with $B$ being the $2^{\infty}$ UHF algebra,
we can require that there be an isomorphism
$\ph \colon A \otimes B \to A.$
Let $\af$ be any automorphism of $A$ which leaves all
tracial states invariant.
Let $\bt$ be an automorphism of $B$ which generates an action
with the Rokhlin property (Proposition~\ref{ExistRP}).
Then $\gm = \ph \circ (\af \otimes \bt) \circ \ph^{-1}$
generates an action of $\Z$ on $A$ with the Rokhlin property,
by Proposition~\ref{TensWithRP}.
Clearly all tracial states on $A$ are $\gm$-invariant,
so Proposition~\ref{SRPImpTRP} shows that $\gm$
generates an action with the \tRp.
Therefore Theorem~\ref{chp3:thm1}, Theorem~\ref{chp4:thm1},
and Theorem~\ref{K0OrdDetByTraces}
show that the crossed product $C^* (\Z, A, \gm)$
has real rank zero and stable rank one,
and that the order on \pj s over this algebra is determined by traces.

We claim that the crossed product $C^* (\Z, A, \gm)$
does not have tracial rank zero.
Using the notation before Definition~3.1 of~\cite{Bn},
we note that the proof of Theorem~7.20 of~\cite{Bn}
gives a tracial state $\ta_0 \in T (A) \setminus T (A)_{\mathrm{IM}}.$
Since $\ta_0$ is assumed to be invariant under $\gm,$
it extends to a tracial state $\ta$ on $C^* (\Z, A, \gm).$
Using the equivalence
of Conditions~(1) and~(4) in Theorem~3.1 of~\cite{Bn},
it follows that
$\ta \in T (C^* (\Z, A, \gm)) \setminus T (C^* (\Z, A, \gm))_{\mathrm{IM}}.$
So $C^* (\Z, A, \gm)$ fails to have tracial rank zero
for the same reason that $A$ does.

We have not determined whether $C^* (\Z, A, \gm)$ is quasidiagonal,
but it seems reasonable to hope that one can use
the \tRp\  to show that it is.
\end{exa}

\begin{exa}\label{NonTAF3}
Let $A$ be a simple separable \ca\  constructed as in
Theorem~7.23 of~\cite{Bn}.
This algebra has real rank zero and stable rank one,
satisfies the local approximation property of Popa and
is hence quasidiagonal (as in Example~\ref{NonTAF2}),
and has a unique tracial state, but is not exact.
It also does not have tracial rank zero.
We show below that $A$ may be chosen such that in addition
the order on \pj s over $A$ is determined by traces.

Let $B$ be the $2^{\infty}$ UHF algebra.
It is easy to see that all the properties given above for $A$
carry over to $A \otimes B.$
(To see that $A \otimes B$ does not have tracial rank zero,
observe that the last paragraph of the proof of Theorem~7.23 of~\cite{Bn}
applies just as well to $A \otimes B$ as to $A.$)
Let $\af$ be any automorphism of $A.$
Let $\bt$ be an automorphism of $B$ which generates an action
with the Rokhlin property (Proposition~\ref{ExistRP}).
Then $\gm = \af \otimes \bt$
generates an action of $\Z$ on $A \otimes B$ with the Rokhlin property,
by Proposition~\ref{TensWithRP}.
Since $A \otimes B$ has a unique tracial state,
it follows from Proposition~\ref{SRPImpTRP} that
$\gm$ generates an action with the \tRp.
Theorem~\ref{chp3:thm1}, Theorem~\ref{chp4:thm1},
and Theorem~\ref{K0OrdDetByTraces} now show that the crossed product
$C^* (\Z, \, A \otimes B, \, \gm)$
has real rank zero, stable rank one, and a unique tracial state,
and that the order on \pj s over this algebra is determined by traces.

We now show how to arrange that
the order on \pj s over $A$ is determined by traces.
This is done by adding one more condition to Conditions~(1) through~(5)
at the beginning of the proof of Theorem~7.23 of~\cite{Bn}.
In addition to the dense sequence $\left\{ \rsz{a_k^{(m)}} \right\}_{k \in \N}$
in $A_m$ used there,
we let $\left\{ \rsz{p_k^{(m)}} \right\}_{k \in \N}$ be a countable
set of \pj s in $\bigcup_{l = 1}^{\infty} M_l \otimes A_m$
such that every \pj\  in $\bigcup_{l = 1}^{\infty} M_l \otimes A_m$
is \mvnt\  to some $p_k^{(m)}.$
Then we require, in addition to Conditions~(1) through~(5),
the existence of a finite set ${\mathcal{P}}_n \subset A_n$ such that
whenever $p$ and $q$ are \pj s in
\[
\ts{
\sm_{n, 0} \left(
  \left\{ \rsz{p_k^{(0)}} \right\}_{1 \leq k \leq n} \right)
\cup
\sm_{n, 1} \left(
  \left\{ \rsz{p_k^{(1)}} \right\}_{1 \leq k \leq n} \right)
\cup \cdots \cup
\sm_{n, \, n - 1}
   \left( \left\{ \rsz{p_k^{(n - 1)}} \right\}_{1 \leq k \leq n} \right) }
\]
such that $\ta (p) < \ta (q)$ for all tracial states $\ta$ on $A_n,$
then there is $s \in \bigcup_{l = 1}^{\infty} M_l \otimes A_n,$
all of whose matrix entries are in ${\mathcal{P}}_n,$
such that $s s^* = p$ and $s^* s \leq q.$
To see that this can be done, at the step in the proof where the sets
${\mathcal{S}}$ and ${\mathcal{I}}$ are chosen,
we observe that the order on \pj s over
$\prod_{j \in \N} M_{k (n (j))} (\C)$ is determined by traces,
choose ${\mathcal{P}}$ accordingly, and include ${\mathcal{P}}$
along with ${\mathcal{U}},$ ${\mathcal{S}},$ and ${\mathcal{I}}$
when generating the next \ca.
\end{exa}

\end{document}